\documentclass[11pt,reqno]{amsart}
\usepackage{amsfonts,amsmath,amsthm,mathtools,bbm,cool}

\usepackage[authoryear,sort]{natbib}
\usepackage{dsfont}
\usepackage{xcolor}
\usepackage{graphicx}
\usepackage{bm}
\usepackage[margin=1in]{geometry}
\usepackage[colorlinks=true,citecolor=magenta,linkcolor=blue]{hyperref}
\usepackage{underscore}

\numberwithin{equation}{section}

\newtheorem{remark}{Remark}[section]
\theoremstyle{plain}

\allowdisplaybreaks

\newcommand{\be}{\begin{equation}}
	\newcommand{\ee}{\end{equation}}

\def\ss{{**}}

\newcommand{\R}{\mathbb{R}}

\newcommand{\Bf}{\mathbf{f}}
\newcommand{\bQ}{\mathbf{Q}}
\newcommand{\bR}{\mathbf{R}}
\newcommand{\CC}{\mathcal{C}}
\newcommand{\fc}{\widehat f}

\newcommand{\fq}{|\widehat{f}|^2}
\newcommand{\dxixi}{\frac{\partial^2}{\partial \xi^2}}
\newcommand{\dxi}{\frac{\partial}{\partial \xi}}
\newcommand{\dt}{\frac{\partial}{\partial t}}

\newenvironment{equations}{\equation\aligned}{\endaligned\endequation}
\newcommand{\fer}[1]{(\ref{#1})}

\begin{document}

\author[G. Martalò]{Giorgio Martalò}
\address{Giorgio Martalò \hfill\break
	Department of Mathematics ``F. Casorati'' \hfill\break
	University of Pavia \hfill\break
	Via Ferrata 5, 27100 Pavia, Italy \vspace*{2mm}}
\email{giorgio.martalo@unipv.it \vspace*{5mm}}

\author[G. Toscani]{Giuseppe Toscani}
\address{Giuseppe Toscani \hfill\break
	Department of Mathematics ``F. Casorati'' \hfill\break
	University of Pavia \hfill\break
	Via Ferrata 5, 27100 Pavia, Italy \hfill\break
\& Institute of Applied Mathematics and Information Technologies (IMATI) \hfill\break National Research Council (CNR), Pavia, Italy \vspace*{2mm}}
\email{giuseppe.toscani@unipv.it \vspace*{5mm}}

\author[M. Zanella]{Mattia Zanella}
\address{Mattia Zanella \hfill\break
	Department of Mathematics ``F. Casorati'' \hfill\break
	University of Pavia \hfill\break
	Via Ferrata 5, 27100 Pavia, Italy \vspace*{2mm}}
\email{mattia.zanella@unipv.it \vspace*{5mm}}

		\title[Individual-Based Foundation of SIR-type epidemic models]{Individual-Based Foundation of SIR-Type Epidemic Models: mean-field limit and large time behaviour}
	
	\date{}
	\maketitle

	\begin{abstract}
We introduce a kinetic framework for modeling the time evolution of the statistical distributions of the population densities in the three compartments of susceptible, infectious, and recovered individuals, under epidemic spreading driven by susceptible–infectious interactions. The model is based on a system of Boltzmann-type equations describing binary interactions between susceptible and infectious individuals, supplemented with linear redistribution operators that account for recovery and reinfection dynamics. The mean values of the kinetic system recover a SIR-type model with reinfection, where the macroscopic parameters are explicitly derived from the underlying microscopic interaction rules. In the grazing collision regime, the Boltzmann system can be approximated by a system of coupled Fokker-Planck equations. This limit allows for a more tractable analysis of the dynamics, including the large-time behavior of the population densities. In this context, we rigorously prove the convergence to equilibrium of the resulting mean-field system in a suitable Sobolev space by means of the so-called energy distance. The analysis reveals the dissipative structure of the dynamics and the role of the interaction terms in driving the system toward a stable equilibrium configuration. These results provide a multi-scale perspective connecting kinetic theory with classical epidemic models.
	\end{abstract}
	
	\smallskip
	
	\noindent{\bf Keywords:} Systems of kinetic equations, Many-agent systems, Mathematical epidemiology
	
	\smallskip
	
	\noindent{\bf MSC Classification:} 35Q84, 35Q20, 35Q82, 82C40, 92D30

	\maketitle

\section{Introduction}\label{sec:intro}

Building upon earlier efforts to understand disease dynamics, almost one century ago Kermack and McKendrick published on this journal the seminal paper \cite{kermack1927} that significantly advanced the study of epidemic spread. In their model, the population is divided into three compartments: susceptible (S), infectious (I), and recovered (R), with their sizes evolving over time. Reinfection dynamics can be captured by suitably modifying classical SIR-type models. In particular, the resulting model takes the form
\begin{equations}\label{SIR-or}
\frac{dS(t)}{dt}&= -\beta\frac{I(t) S(t)}N +\alpha R(t), \qquad S(0) = S_0 \ge 0,\\
\frac{dI(t)}{dt}&= \beta\frac{I(t) S(t)}N - \gamma I(t), \qquad I(0) = I_0 \ge 0, \\
\frac{dR(t)}{dt}&=  \gamma I(t) -\alpha R(t), \qquad I(0) = R_0 \ge 0,
\end{equations}
where $S(t), I(t)$, and $R(t)$ are the numbers of individuals in these compartments, satisfying the relation $S(t)+I(t)+R(t) = N$ for all $t\ge0$. In \fer{SIR-or}   $\beta>0$ measures the average infection rate of the population and $\gamma>0$ is the recovery rate in the compartment of infectious agents, being $\alpha>0$  the loss of immunity rate. System \fer{SIR-or} belongs to the class of epidemiological models that rely on a compartmentalization of the population, denoted by $\mathcal C$, where individuals are divided into epidemiologically relevant categories. Nevertheless,  other compartments may contribute to the transmission of the disease \cite{hethcote2000}. Hence, alternative modeling approaches are possible by adopting different forms of compartmentalization, depending on the characteristics of the disease or the specific aspects of interest.

These models serve as important conceptual tools for understanding disease propagation in a uniformly mixed population, with transition rates between compartments typically postulated based on empirical evidence or heuristic reasoning. However, since they are not derived from first principles, they come with inherent limitations. In practice, mixing patterns are influenced by numerous factors, such as population density and social behaviour, which lead to significant heterogeneity in the number of contacts individuals have.

For these reasons, while maintaining the division into  compartments, mathematical research started to move toward a more precise description of the temporal evolution of variations in the classes of individuals in the presence of multiple time-scales. In particular, methods from classical kinetic theory appear to be well suited to describe the evolution of individuals across compartments in terms of their statistical distribution with respect to additional variables that quantify certain population features, such as possible in-host heterogeneities, agent-based number of contacts, or an evolving opinion variable \cite{dimarco2021kinetic,albi2022kinetic,zanella2023,bertaglia2023,dellamarca23,lorenzi24,BisiLorenzani,bertaglia24,bernardi2025heterogeneously}. This framework allows to model how these and other distributions evolve in the presence of an epidemic also in a space dependent setting, see e.g. \cite{ha22,Difrancesco25}. In this context, few results are today available regarding the passage from microscopic trajectories to averaged quantities where,  in different communities, huge research efforts have been devoted to derive mean-field models. These equations encapsulate individual-based dynamics in the limit of a large number of agents and can be obtained with rigorous limiting procedures, see \cite{bolley12,carrillo10,degond08,ha09,motsch14}. Indeed, it is a well known fact that the main disadvantage of the microscopic approach to capture the emerging behavior of interacting systems relies on the so-called curse of dimensionality. Despite the loss of microscopic information, these aggregate descriptions enable the study of emergent phenomena by linking individual-based dynamics with observable macroscopic quantities and phase-transitions, see \cite{during09,during15,MR4251319,MR3067586,preziosi}. 

In what follows, we deal with a mean-field description of the compartmental epidemic model in the form 
\be\label{vec-LV}
\frac{\partial\Bf(x,t)}{\partial t} =  \bQ(\Bf(x,t)) + \bR(\Bf(x,t)), \qquad  \Bf(x,t=0) = \Bf_0(x),
\ee	
where $x \in \mathbb R_+$.  In \fer{vec-LV} the vector $\bQ$ is the $n$-dimensional vector whose components $Q_j( \Bf)$, $j \in \mathcal C$ are suitable (generally nonlinear) kinetic  operators which describe the variation of the densities $f_j(x,t)$ of individuals  in the $j$-th compartment due  to interactions with individuals
in other compartments, while  $\bR(\Bf(x,t))$ is the $n$-dimensional vector whose components $R_j( f_j)$, $j \in \mathcal C$ are linear kinetic redistribution operators which quantify the variation of the $f_j$ density due to passage towards and/or arrival from other compartments \cite{bisi2009,BonMenTosZan,TosZan}. System \fer{vec-LV} is built to be closely related to classical epidemiological models, which represent nothing but the evolution of the mean values of the densities $f_j(x,t)$, $j \in \mathcal C$. 

Motivated by the analogies between epidemic spreading and predator–prey dynamics, and building on the recent results in \cite{BonMenTosZan,TosZan}, we consider here the kinetic formulation \fer{vec-LV} of the classical SIR model with reinfection \fer{SIR-or}, in the setting where the variable $x\in \mathbb R_+$  represents the spatial density of the population, which is related to the incidence of the disease, see e.g. \cite{DeJD,HU2013125,riley07}. 

Our main result concerns the time evolution of population densities in two cases: the simpler scenario without reinfection ($\alpha = 0$) and the general case with reinfection ($\alpha>0$). In the absence of reinfection, we show that starting from any initial configuration with at least two bounded moments, the epidemic dynamic leads to an explicit modification of the population densities in the three compartments, which retain a memory of the initial state at all subsequent times. In contrast, when reinfection is present, this memory is asymptotically lost, and the system's densities converge to computable steady states determined by the parameters governing the epidemic.

At variance with classical entropy methods for single Fokker-Planck equations, see e.g. \cite{Furioli,Toscani99,Toscani:1999aa}, we cannot obtain a sharp equilibration rate due to the intrinsic coupling of the Fokker-Planck system of equations. To this end, we will provide explicit upper bounds for the equilibration rates for the considered system of coupled Fokker-Planck equations in a homogeneous Sobolev space $\dot{H}_{-p}$ with $1/2<p\le 3/2$. This norm can be equivalently expressed in the physical space through the so-called energy distance, of common use in statistics \cite{szekely2013}. 
These results can be used to assess the impact of population heterogeneity on the final steady state of an epidemic and to calibrate efficiently the external interventions with a more accurate representation of the population. Such interventions effectively alter the shape of the initial population densities, thereby influencing the long-term dynamics of the system.

In more detail, the paper is organized as follows. Section \ref{sec:kinetic} is devoted to a precise description of the elementary interactions which enable to build up the kinetic model of epidemic spreading, which is dealt with in Section \ref{sec:Boltzmann}.  The main feature of this model, which relies on the evolution of the population's distributions in terms of their territorial densities, is its connection with the classical SIR model \fer{SIR-or}, which appears when evaluating the mean values of the statistical distribution of the involved populations. Then, Section \ref{sec:FP} introduces a simplification of the kinetic model, in the form of a system of Fokker--Planck type equations with variable in time coefficients of diffusion and drift. The large--time evolution of the solution to the system of Fokker–Planck type equations is subsequently studied in Section \ref{sec:large}, where it is shown that the alternative cases without and with re-infection lead to a completely different behavior. Last, numerical simulations illustrating the effectiveness of the approach are presented in Section \ref{sec:numerics}.

\section{A kinetic description of epidemic spreading}\label{sec:kinetic}

We consider a system composed of three populations, which represent the possible states of individuals with respect to the  epidemic spreading: susceptible $S$, infectious $I$ and recovered $R$. In short, we will indicate by $\mathcal C=\{S,I,R\}$ the set of three populations. We assume that each population is characterized at a statistical level by its distribution with respect to a relevant variable, denoted by $x\in\mathbb{R}_+$, which quantifies the \emph{spatial density} of individuals.  This variable, owing to the fact that big cities  are characterized by a high density of people in a relatively small area,  can be easily related to the size distribution of cities, which is generally well described by a probability density with polynomial tails (cfr. \cite{gualandi2019} and the references therein). This choice is clearly motivated by the fact that a typical measure to contain the epidemic spreading is to resort to lockdown policies, which correspond to rarefying the average number of contacts between agents.

Let us denote by  $f_J(x, t)$, $J \in \mathcal C$ the distribution functions, which are such that $f_J(x,t)\,dx$ indicates the fraction of individuals of population $J$ belonging to the interval $(x, x +dx)$, with $x\ge0$, at time $t$. Then, having in mind that we want to describe the epidemic spreading in terms of the population density,  we assume that the evolution over time $t\ge 0$ of the statistical distribution of the  populations of susceptible and infectious individuals is fully dependent  of its  elementary interactions with individuals of the other population. 

Without loss of generality, we can assume that the initial distribution functions (at time $t=0$) are probability density functions, i.e., for $J \in \CC$
\begin{equation}
    \int_{\mathbb{R}_+}f_J(x,0)\,dx\,=\,1\,.
    \label{intial_condition}
\end{equation}
In this way, assuming that the distributions functions are mass preserving, we can consistently define, for any time $t \ge 0$,  the principal moments of order $k\ge0$
\begin{equation}
\label{eq:moment}
    m_J^{(k)}(t)\,:=\int_{\mathbb{R}_+}x^kf_J(x,t)\,dx.
\end{equation}
The case $k= 1$ defines the (variable in time) mean value of the group $J$. In this case, to simplify notations , we define, for $J\in \CC$, 
\be\label{eq:mean}
m_J(t) = m_J^{(1)}(t) := \int_{\mathbb{R}_+}xf_J(x,t)\,dx\,.
\ee
The case $k = 2$ is useful to quantify, through the notion of variance
\be\label{eq:variance}
V_J(t) \, := \int_{\mathbb{R}_+}\left( x - m_J(t)\right)^2\,f_J(x,t)\,dx = m_J^{(2)}(t) - [m_J(t)]^2 ,
\ee
the spreading of the populations around their mean values.  It is important to remark that identity \fer{eq:variance} allows to express moments of order two (cf. \fer{eq:moment} for $j=2$) in terms of the variance and of moments of order one
\be\label{eq:energy}
m_J^{(2)}(t)\,  = V_J(t) +  [m_J(t)]^2.
\ee

\subsection{The elementary interactions}\label{sec:elementary}

In the epidemic spreading, the relevant elementary interactions involve contacts between susceptible and infectious individuals. The outcome of these interactions need to be designed in agreement with the expected behavior of the  system under study. First, susceptible--infectious (S-I) interactions should lead to a decrease in the number of susceptible individuals and, conversely, to an increase in the number of infectious. Second, since  it is not certain that the result of the interaction will lead to infection of the susceptible,  the outcome should also be subject to random effects. This type of interactions is  analogous to a binary interaction between the populations of preys (here the susceptible individuals) and predators (the infectious individuals) in a Lotka--Volterra type kinetic description, a situation which has recently been considered in \cite{TosZan,BonMenTosZan}.  Following the natural  assumptions in \cite{TosZan,BonMenTosZan}, we define the instantaneous microscopic transition $x \mapsto x'$ and $x_*\mapsto x_*'$ between a susceptible individual characterized by a  value $x$ and an infectious individual characterized by a value $x_*$  in the form
\begin{equation} \label{cross_interactions}
    \begin{split}
        x' &= x - \Phi(x_*)x +  \eta_{{S}}(x_*)\, x,
        \\[2mm]
        x_*' &= x_* + \Phi(x)x_* +  \eta_{{I}}(x)\, x_*\,;
    \end{split}
\end{equation}
in \fer{cross_interactions}, the function $\Phi(y)$, which quantifies deterministic variations within the S-I interactions, is assumed of the following form  
\begin{equation} \label{eq:deterministic effects}
        \Phi(x_*) = \beta \frac{x_*}{1+x_*}.
\end{equation}
In the first equation of \fer{cross_interactions}, this increasing function models  the infection rate. This function assumes the value zero in absence of infectious individuals, and, in agreement with classical approaches to epidemic modeling \cite{capasso1978, liu1986},  takes into account saturation phenomena for large numbers of infectious individuals. Clearly the same function quantifies, in the second equation of \fer{cross_interactions}, the increase in the number of infectious, and coherently it is assumed proportional to the amount of susceptible agents having contacts with them. Here, $\beta \in (0,1)$ is a positive constant measuring the maximal infection rate.  Note in particular that  the infection rate remains always smaller than one.  Such type  of  interactions are also in agreement with the biological literature on similar subjects, where this behavior is know under the name of Holling-type growth \cite{holling1959components, holling1966functional, beddington1975mutual, kuno1987principles, abrams2000nature}. 

Relations \fer{cross_interactions} also involve stochastic variations, which are modeled by  two independent random variables $\eta_{{S}}$ and $\eta_{{I}}$, having zero mean and variable bounded variances of the form
\begin{equation} \label{eq:random effects}
        \langle\eta_{{S}}^2(x_*)\rangle = \sigma_S\frac{x_*}{1+x_*}, \quad
        \langle\eta_{{I}}^2(x)\rangle = \sigma_I\frac{x}{1+x}.
\end{equation}
In \fer{eq:random effects} we have denoted with $\langle \cdot \rangle$ the integration with respect to the probability densities of $\eta_S$ and $\eta_I$,  with $\sigma_S, \sigma_I$ some constant positive parameters. These random effects  encode   
the possible unpredictable variations of the elementary process of epidemic spreading. According to this assumption,  the distribution of the randomness around the mean behavior of a population is assumed proportional to the amount of the other populations. Also in this case, the choice  of a dependence from the densities like \fer{eq:deterministic effects} for the variances of $\eta_S$ and $\eta_I$ translates the fact that the random effects can not exceed (in mean) a saturation value prescribed by the constants $\sigma_S$ and $\sigma_I$. 

\begin{remark}\label{rem:positivity}
 
It is essential that the microscopic transitions \eqref{cross_interactions} must guarantee the positivity of post-interactions sizes $x^\prime$ and $x_*^\prime$. 
Since $\Phi(x_*) \le \beta$, positivity of the post-interaction value $x'$ requires that the negative values assumed by the random variable $\eta_1(x_*)$ are greater than $-(1-\beta)$.
Likewise, since $\Phi(x) \ge 0$, positivity of the post-interaction value $x_*'$ requires that the negative values assumed by the random variable $\eta_2(x)$ are greater than $-1$. 
\end{remark}


In addition to the elementary interaction responsible of the epidemic spreading,  other important processes need to be taken into account.  A leading aspect to be considered in the modeling of an epidemic spreading relies in the recovery dynamics.  At the elementary level, the dynamics of the recovery of the infected population and the loss of immunity of recovered agents, that describe from one side the statistical effect of drug treatment, and from the other side the statistical description of  new disease mutations  and/or the loss of immunity over time, can be described by resorting to background-type interactions. 

The recovery process can be described in the following way: given $x_* \in X_I \sim f_I(x,t)$ we consider the linear interaction
\begin{equation}
    \label{eq:redistribution1}
    \begin{split}
        x_*^{\prime\prime}& = x_*+\gamma(y-\theta x_*), 
        \end{split}
\end{equation}
where $y \sim h(y,t)$ represent a background distribution encapsulating the availability of medical treatments. This distribution is such that
\[
\int_{\mathbb R_+}yh(y,t)dy = (\theta-1)m_I(t),
\]
where the constant $\theta >1$. In \fer{eq:redistribution1} the positivity of post-interaction state $x_*^{\prime\prime}$ follows under the condition $\gamma\theta<1$. Similarly, the gain of the recovered agents are such that, given $x_\ss \in X_R \sim f_R(x,t)$ and $x_* \in X_I \sim f_I(x,t)$
\begin{equation}
\label{eq:redistribution2}
x_\ss^{\prime} = x_\ss + \gamma x_*.
\end{equation}
Likewise,  the loss of immunity of recovered agents is described by a linear interaction that takes into account the randomness of the process. Given $x_\ss \in X_R \sim f_R(x,t)$, this interaction takes the form
\begin{equation}
\label{eq:redistribution3}
x_\ss^{\prime\prime} = x_\ss - \alpha x_\ss + x_\ss \eta_R, \qquad 0\le\alpha<1,
\end{equation}
{being $\eta_R$ a random variable with that $\left\langle \eta_R \right\rangle = 0$ and $\left\langle \eta_R^2 \right\rangle = \sigma_R(\alpha)$, where $\sigma_R(\alpha=0)= 0$. }

Last, the passage to the population of susceptible due to the loss of immunity is obtained by considering, for $x \in X_S \sim  f_S(x,t)$ and $x_\ss \in X_R \sim f_R(x,t)$,  the transition
\begin{equation}
\label{eq:redistribution4}
x^{\prime\prime} = x + \alpha x_\ss.
\end{equation}

\begin{remark}\label{rem:loss}
The previous description of the process of loss of immunity includes both the case of re-infection, and in this case $\alpha >0$, and the case without re-infection, which corresponds to the choice $\alpha =0$. Clearly, in this last case the interaction pair 
\fer{eq:redistribution3}, \fer{eq:redistribution4} does not produce any variation in the populations of recovered and susceptible individuals. 
\end{remark}

\section{The kinetic model}\label{sec:Boltzmann}
Once the details of the elementary interactions have been characterised,  the evolution of the distribution functions $f_J(x,t)$, $J\in \CC$ is governed by a  Boltzmann-type system of equations. Following the approach in \cite{MR1307620,ParTos-2013}, we get
\begin{equations}\label{eq:kine}
    \dfrac{\partial f_S(x,t)}{\partial t}&=\mathcal{Q}_{SI}(f_S,f_I)(x,t)+\mathcal{R}^\alpha(f_S,f_R)(x,t),\\
    \dfrac{\partial f_I(x,t)}{\partial t}&=\mathcal{Q}_{IS}(f_I,f_S)(x,t)+ \mathcal{L}^{-\gamma}(f_I)(x,t), \\
    \dfrac{\partial f_R(x,t)}{\partial t}&=\mathcal{R}^{\gamma}(f_I,f_R)(x,t)+\mathcal{L}^{-\alpha}(f_R)(x,t).
\end{equations}
In system \fer{eq:kine} the operators $\mathcal{Q}_{SI}, \mathcal{Q}_{IS}$ are bilinear Boltzmann operators based on the cross interactions \eqref{cross_interactions} and read
\[
\begin{split}
\mathcal{Q}_{SI}(x,t) = \left\langle\int_{\mathbb R_+} \kappa(x_*) \left(\dfrac{1}{{}^\prime J_{SI}}f_S({}^\prime x,t)f_I({}^\prime x_*,t) - f_S(x,t)f_I(x_*,t)\right)dx_*\right\rangle\,, \\
\mathcal{Q}_{IS}(x,t) = \left\langle\int_{\mathbb R_+} \kappa(x_*) \left(\dfrac{1}{{}^\prime J_{IS}}f_I({}^\prime x,t)f_S({}^\prime x_*,t) - f_I(x,t)f_S(x_*,t)\right)dx_*\right\rangle. 
\end{split}\]  
where $\left\langle \cdot \right\rangle$ is the expectation with respect to the random variables $\eta_S,\eta_I$ and ${}^\prime J_{SI}, {}^\prime J_{IS}$ are the determinants of the Jacobian matrices characterising the transition $({}^\prime x,^{\prime}x_*) \to (x,x_*)$ as in in \eqref{cross_interactions}, while the operator $\mathcal{R}^\alpha(f_S,f_R)$ represents the re-infection due to the microscopic interaction dynamics \eqref{eq:redistribution4} and describing the transition $(^{\prime\prime}x,x_\ss) \to (x,x_\ss)$. In particular, this operator reads
\[
\begin{split}
\mathcal R^\alpha(f_S,f_R)(x,t) =  \int_{\mathbb R_+} \left(\dfrac{1}{{}^{\prime\prime} J_{SR}} f_S({}^{\prime\prime} x,t)f_R( x_{**},t) - f_S(x,t)f_R(x_{**},t)\right)dx_{**}, 
\end{split}
\]
being ${}^{\prime\prime}J_{SR}$ the determinant of the Jacobian associated to the transformation ${}^{\prime\prime}x \to x$, whereas the operator $\mathcal{R}^{\gamma}(f_I,f_R)$ encapsulates the recovery due to the microscopic dynamics \eqref{eq:redistribution2} for the transition $(x_*,^{\prime\prime}x_\ss) \to (x_*,x_\ss)$ and reads
\[
\begin{split}
\mathcal R^\gamma(f_I,f_R)(x,t) =  \int_{\mathbb R_+} \left(\dfrac{1}{{}^{\prime\prime} J_{IR}} f_I( x_*,t)f_R( {}^{\prime\prime}x,t) - f_I(x_*,t)f_R(x,t)\right)dx_*. 
\end{split}
\]
The linear operator $\mathcal L^{-\gamma}(f_I)$ represents the loss of infection following the microscopic transition $^{\prime\prime}x_*\to x_*$ defined in \eqref{eq:redistribution1} and reads
\[
\begin{split}
\mathcal L^{-\gamma}(f_I)(x,t) =\int_{\mathbb R_+^2} \dfrac{1}{{}^{\prime\prime} J_{I}} f_I({}^{\prime\prime} x,t)h(y,t)dydx - f_I(x,t),
\end{split}
\]
being ${}^{\prime\prime} J_I$ the determinant of the Jacobian associated to \eqref{eq:redistribution1}. Finally the operator $\mathcal L^{-\alpha}(f_R)$ encapsulates the loss of immunity due to the microscopic transition $^{\prime\prime}x_\ss \to x_\ss$ defined in \eqref{eq:redistribution3} and reads
\[
\begin{split}
\mathcal L^{-\alpha}(f_R)(x,t) =\left\langle\int_{\mathbb R_+} \dfrac{1}{{}^{\prime\prime} J_{R}} f_R({}^{\prime\prime} x,t)dx\right\rangle - f_R(x,t),
\end{split}
\]
being ${}^{\prime\prime} J_R$ the determinant of the Jacobian associated to \eqref{eq:redistribution3}.  To investigate the action of system \fer{eq:kine} on the observable macroscopic quantities, it is useful to consider the weak form of the Boltzmann equations.  It corresponds to write, for any given set of smooth functions  $\varphi_J(x)$, $J \in \mathcal{C}$,  the system 
of equations

\begin{equation}
\label{weak_S}
\begin{split}
   & \dfrac{d}{dt}\int_{\mathbb{R}_+}\varphi_S(x)f_S(x,t)dx=  \underbrace{\int_{\mathbb{R}_+^2}\kappa(x_*)\langle\varphi_S(x^{\prime})-\varphi_S(x)\rangle f_S(x,t)f_I(x_*,t)dx\, dx_*}_{\mathcal{I}_{S}} \\
    &\quad+\underbrace{\int_{\mathbb{R}_+^2}(\varphi_S(x^{\prime\prime})-\varphi_S(x))f_S(x,t)f_R(x_\ss,t) dx\,dx_\ss}_{\mathcal{J}_S}\,,
\end{split}
\end{equation}
    
\begin{equation}
\label{weak_I}
\begin{split}
&\dfrac{d}{dt}\int_{\mathbb{R}_+}\varphi_I(x_*)f_I(x_*,t)dx_*=\underbrace{\int_{\mathbb{R}_+^2}\kappa(x)\langle\varphi_I(x_*^\prime)-\varphi_I(x_*)\rangle f_S(x,t)f_I(x_*,t)dxdx_*}_{\mathcal I_{I}}\\
&\quad+\underbrace{\int_{\mathbb{R}_+^2}(\varphi_I(x_*^{\prime\prime})-\varphi_I(x_*))f_I(x_*,t)h(y,t) dx_*dy }_{\mathcal J_I}\,,
\end{split}
\end{equation}
    
\begin{equation}
\label{weak_R}
\begin{split}
   & \dfrac{d}{dt}\int_{\mathbb{R}_+}\varphi_R(x_\ss)f_R(x_\ss,t)dx_\ss=\underbrace{\int_{\mathbb{R}_+^2}\left(\varphi_R(x_\ss^{\prime})-\varphi_R(x_\ss)\right) f_I(x_*,t)f_R(x_\ss,t)dx_*\, dx_\ss}_{\mathcal I_R} \\
    &+\underbrace{\int_{\mathbb{R}_+}\left\langle\varphi_R(x_\ss^{\prime\prime})-\varphi_R(x_\ss)\right\rangle f_R(x_\ss,t)dx_\ss}_{\mathcal J_R}\,.
\end{split}
\end{equation}
In \fer{weak_S} and \fer{weak_I} the notation $\langle \cdot\rangle$ denotes the mean value of the interaction integral with respect to the random quantities $\eta_i$, $i =1,2$ present in the elementary variations \fer{cross_interactions}.
Note moreover that, in accord with classical kinetic theory,  the  interaction integrals in  \fer{weak_S} and \fer{weak_I}  depend on the interaction frequency $\kappa$, which for susceptible individuals is consistently assumed proportional to the local density of the infectious ones, and analogously for the infectious individuals, with 
\be
\label{eq:kappa}
\kappa(x_*)=1+x_*; \quad \kappa(x)=1+x .
\ee
\subsection{Evolution of the main moments }\label{sec:moments}

The weak form of system  \fer{eq:kine} allows to compute the principal moments. It is immediate to conclude that the choice $\varphi_S(x)=\varphi_I(x)=\varphi_R(x)=1$,  leads to mass conservation. Therefore, owing to \eqref{intial_condition}, it follows that
\begin{equation}
    \int_{\mathbb{R}_+}f_S(x,t)dx=\int_{\mathbb{R}_+}f_I(x,t)dx=\int_{\mathbb{R}_+}f_R(x,t)dx=1
    \label{density_at_any_time}
\end{equation}
for any time $t\in\mathbb{R}_+$. In view of this normalisation,  $f_J$, $J \in \mathcal{C}$, are probability densities. The present framework departs from earlier kinetic formulations since its structure establishes a mathematically consistent link between the underlying agent-based dynamics and the emergent collective behaviour described through the moments of these kinetic densities. The choice $\varphi_J(x)=x$, $J \in \mathcal{C}$ allows to compute the evolution of the mean values \fer{eq:mean}, i.e. $m_S(t),m_I(t),m_R(t)$.  The balance between the interaction frequency   \eqref{eq:kappa} and  the elementary  interactions defined in \eqref{cross_interactions} and \eqref{eq:redistribution4}  gives
\begin{equation}
\begin{split}
    \dfrac{d}{dt}m_S(t) &=\int_{\mathbb{R}_+^2}(1+x_*)\langle x^{\prime}-x\rangle f_S(x,t)f_I(x_*,t)dxdx_* + \\
    &\qquad \int_{\mathbb{R}_+^2}(x^{\prime\prime}-x)f_S(x,t)f_R(x_\ss,t)dx\,dx_\ss \\
    &=-\beta \, m_S(t) m_I(t) + \alpha \, m_R(t)\,.
\end{split}
\end{equation}
Likewise, if $\varphi_I(x)=x$, we obtain
\begin{equation}
\begin{split}
\dfrac{d}{dt}m_I(t) &=\int_{\mathbb{R}_+^2}(1+x)\langle x_*^{\prime}-x_*\rangle f_S(x,t)f_I(x_*,t)dx\, dx_* +\\
&\qquad  \int_{\mathbb{R}_+}(x_*^{\prime\prime}-x_*)f_I(x_*,t)h(y,t)dx_*dy\\
    & = \beta \, m_S(t) m_I (t) - \gamma \, m_I(t)\,.
    \end{split}
\end{equation}
Setting finally $\varphi_R(x) = x$, we obtain
\begin{equation}
\begin{split}
\dfrac{d}{dt}m_R(t) &=\int_{\mathbb{R}_+^2}(x_\ss^{\prime}-x_\ss)f_I(x,t)f_R(x_\ss,t)dx\, dx_\ss + \\
&\qquad \int_{\mathbb{R}_+}\left\langle x_\ss^{\prime\prime}-x_\ss\right\rangle f_R(x_\ss,t) dx_\ss \\
	&=\gamma \, m_I(t)-\alpha \, m_R(t)\,.
	\end{split}
\end{equation}
Therefore, the evolution of the mean value obeys to the classical SIR model with re-infection \fer{SIR-or}
\begin{equation}
\label{eq:mSIR}
\begin{split}
\dot m_S(t) &= -\beta m_S(t)m_I(t)+ \alpha m_R(t)\,, \\
\dot m_I(t) &= \beta m_S(t)m_I(t)- \gamma m_I(t)\,, \\
\dot m_R(t)&= \gamma m_I(t) - \alpha m_R(t)\,,
\end{split}
\end{equation}
whose equilibrium, under the natural conditions  $\gamma < \beta$, and $\alpha>0$,  is given by \cite{iannellipugliese} 
\begin{equation}
\label{eq:minfty}
\left(m_S^{\infty},m_I^{\infty},m_R^{\infty}\right) = \left(\frac{\gamma}{\beta},\frac{\alpha}{\gamma+\alpha}\left(1-\frac{\gamma}{\beta}\right) , \frac{\gamma}{\gamma+\alpha}\left(1-\frac{\gamma}{\beta}\right)\right)\,.
\end{equation}

\begin{remark}
The previous result shows that SIR type models represent nothing but the evolution of the mean values of the system of probability densities $f_J(x,t)$, $J\in \CC$. These densities track the epidemic process across classes of contact activity, linked to territorial density, and thus contain more information than the aggregated SIR variables.
Clearly, the picture provided by the densities $f_J(x,t)$ is more detailed than the simple evolution furnished by the mean values, since it allows for the description of the evolution of higher moments. 
Earlier kinetic approaches offered a distinct derivation of SIR-type equations, in which fast-scale fluctuations in social heterogeneity govern the evolution of the mass fractions of kinetic densities. The corresponding higher-order moments naturally encode behavioural variability and its influence on transmission dynamics, see e.g. \cite{dimarco2021kinetic,zanella2023}.
\end{remark}

\begin{remark} \label{rem:mean}
For future purposes, it is important to remark that, under the conditions  $\gamma < \beta$, i.e. when the basic reproduction number $\mathcal{R}_0=\beta/\gamma>1$, and $\alpha>0$, the mean values satisfy, for all times $t \ge 0$  the lower bounds $m_J(t) \ge c >0$, $J \in \mathcal{C}$ \cite{hethcote2000}.
\end{remark}

Further moments of  interest are the variances of the densities, which quantify the spreading around the mean values.  To compute the evolution of the variances in each compartment,  defined in \fer{eq:variance}, 
we consider as test function $\varphi_J(x) = (x-m_J(t))^2$, $J \in \mathcal{C}$.  

 Owing to \eqref{weak_S}-\eqref{weak_I}-\eqref{weak_R} and using \fer{eq:energy}, we get for the compartment of susceptible individuals
\begin{equation}
\label{eq:VS_B}
\begin{split}
	\dfrac{d}{dt} V_S(t) =&  \left[V_S(t) + m^2_S(t)\right]\left\{\int_{\mathbb R_+} \beta^2\dfrac{ x_*^2}{1+x_*}  f_I(x_*,t)\,dx_* + \sigma_S m_I (t)\right\}  \\
	&   - 2\beta m_I(t) V_S(t)+\alpha^2 \left[V_R(t) + m^2_R(t)\right].
\end{split}
\end{equation}
Analogous computations give for the infected compartment
\begin{equation}
\label{eq:VI_B}
\begin{split}
	\dfrac{d}{dt} V_I(t) =&  (V_I + m_I^2)\left\{\int_{\mathbb R_+}\beta^2 \dfrac{x^2 }{(1+x)} f_S(x,t)dx + (2\beta + \sigma_I) m_S(t)\right\}   \\
	&+(\gamma^2\theta^2-2\theta\gamma) (V_I  + m_I^2)+ 2\gamma\theta(1-\gamma(\theta-1))m_I^2 + \gamma^2 \int_{\mathbb R_+}y^2h(y,t)dy.
\end{split}
\end{equation}
Finally,  taking into account Remark \ref{rem:loss}, for the recovered compartment we get
\begin{equation}
\label{eq:VR_B}
\dfrac{d}{dt} V_R(t) = \gamma^2 (V_I + m_I^2) - 2\alpha V_R + (\alpha^2+\sigma_R(\alpha))(V_R + m_R^2). 
\end{equation}

We may observe how the evolution of the variances involve the kinetic distributions $f_J(x,t)$, $J \in \mathcal{C}$ making their dynamics, unlike the ones of the mean values, not closed. 

\section{Reduced  complexity  models}\label{sec:FP}

As shown in Section \ref{sec:moments}, while the kinetic description allows for an exhaustive knowledge, at least from a numerical viewpoint,  of the time evolution of the solutions, except for the evolution of the mean value, the evolution of higher moments is not closed. 
A simplified version of the kinetic system can however be obtained by resorting to the so-called \emph{quasi--invariant} scaling, which is obtained by studying the kinetic system in a regime where the interactions  are \emph{grazing} (which corresponds to impose that there is only an order $\varepsilon \ll 1$ modification of the pre--interactions values), and at the same time a scaling of time is needed to maintain the correct variation of the time-evolution of the solution.

Given the scaling parameter $\varepsilon\ll 1$,  we consider a new scaled time variable $\tau=\varepsilon t$, together with the scaled distribution functions $f^\varepsilon_J(x,\tau)=f_J(x,t/\varepsilon)$ ($J\in\mathcal{C}$) and $h^\varepsilon(y,\tau)=h(y,t/\varepsilon)$, $r^\varepsilon(z,\tau)=r(z,t/\varepsilon)$ for populations and background, respectively. At the same time we scale the relevant parameters
\begin{eqnarray}\label{eq:scaling}
    \beta\longrightarrow\varepsilon\beta\,,\qquad \sigma_J\longrightarrow\varepsilon\sigma_J\,,\qquad\alpha\longrightarrow\varepsilon\alpha\,,\qquad\gamma\longrightarrow\varepsilon\gamma. 
\end{eqnarray}
This scaling is such that the variation in interactions are very small. This allows to resort to a third order Taylor expansion 
\begin{eqnarray*}
    \varphi_J(x^\prime)-\varphi_J(x)&=& \dfrac{d}{dx}\varphi_J(x)(x^\prime-x)+\dfrac12\dfrac{d^2}{dx^2}\varphi_J(x)(x^\prime-x)^2+\dfrac16\dfrac{d^3}{dx^3}\varphi_J(\hat{x})(x^\prime-x)^3, 
    \\
    \varphi_J(x^{\prime\prime})-\varphi_J(x)&=& \dfrac{d}{dx}\varphi_J(x)(x^{\prime\prime}-x)+\dfrac12\dfrac{d^2}{dx^2}\varphi_J(\tilde x)(x^{\prime\prime}-x)^2, 
\end{eqnarray*}
in \eqref{weak_S}, \eqref{weak_I} and \eqref{weak_R}, where 
\begin{eqnarray}
&&\hat{x}\in(\min\{x,x^\prime\},\max\{x,x^\prime\})\,,\qquad \tilde{x}\in(\min\{x,x^{\prime\prime}\},\max\{x,x^{\prime\prime}\})\,.
\end{eqnarray}
Using these expansions into the weak form of the kinetic equations gives 
\begin{eqnarray}
    \mathcal{I}_S&=&\int_{\mathbb{R}_+^2}\kappa(x_*)\langle\varphi_S(x^{\prime})-\varphi_S(x)\rangle f_S^\varepsilon(x,\tau)f_I^\varepsilon(x_*,\tau)dx\,dx_*\nonumber\\
    &=&\int_{\mathbb{R}_+^2}\kappa(x_*)\left(-\varepsilon\beta\dfrac{x \, x_*}{1+x_*}\right)\dfrac{d}{dx}\varphi_S(x)f_S^\varepsilon(x,\tau)f_I^\varepsilon(x_*,\tau)dx\, dx_*\nonumber\\
    &+&\dfrac12\int_{\mathbb{R}_+^2}\kappa(x_*)\left[\varepsilon^2\beta^2\dfrac{x^2\,x_*^2}{(1+x_*)^2}+\varepsilon\sigma_S\dfrac{x^{2}x_*}{1+x_*}\right]\times\nonumber\\
    &&\hspace{3cm}\times\dfrac{d^2}{dx^2}\varphi_S(x)f_S^\varepsilon(x,\tau)f_I^\varepsilon(x_*,\tau)dx\,dx_*    +{R}_S\,, \nonumber
\end{eqnarray}
where 
\begin{equations}\nonumber   
 {R}_S(f_S^\varepsilon,f_I^\varepsilon)&=\dfrac16\int_{\mathbb{R}_+^2}\kappa(x_*)\left\langle\left(-\varepsilon\beta\dfrac{x\,x_*}{(1+x_*}+\eta(x_*)x\right)^3\right\rangle\times \\
    &\hspace{3cm}\times\dfrac{d^3}{dx^3}\varphi_S({x})f_S^\varepsilon(x,\tau)f_I^\varepsilon(x_*,\tau)dx\,dx_*\,.
    \end{equations}
    

If we assume that the random variable $\eta_S$ has bounded third order moment, then we can write $\eta(x_*)=\sqrt{\dfrac{\varepsilon\sigma_Sx_*}{1+x_*}}\,\hat{\eta}$, with $\langle\hat{\eta}\rangle=0$ and $\langle\hat{\eta}^2 \rangle=1$. Then, it follows that  $ |{R}_S(f_S^\varepsilon,f_I^\varepsilon)|\approx \varepsilon^3+\varepsilon^2+\varepsilon\sqrt{\varepsilon}$.
Analogously, for the term $\mathcal J_S$ we get
\begin{equations}\nonumber
    \mathcal{J}_S&=&\int_{\mathbb{R}_+^2}(\varphi_S(x^{\prime\prime})-\varphi_S(x))f_S^\varepsilon(x,\tau)f_R^\varepsilon(x_\ss,\tau)dx\,dx_\ss 
    \\
    &=& \varepsilon\int_{\mathbb{R}_+^2}\alpha x_\ss \dfrac{d}{dx}\varphi_S(x)f_S^\varepsilon(x,\tau)f_R^\varepsilon(x_\ss,\tau)dx\,dx_\ss \,.
\end{equations}
Letting $\varepsilon \to 0$, and considering that higher order terms in $\varepsilon \to 0$ vanish,  one shows that the distribution $f_S^\varepsilon$ ia weakly convergent to $f_S$,  solution of   the  equation 
\begin{eqnarray}
    &&\dfrac{d}{d\tau}\int_{\mathbb{R}_+}\varphi_S(x)f_S(x,\tau)dx=-\beta m_I(t)\int_{\mathbb{R}_+}x\dfrac{d}{dx}\varphi_S(x)f_S(x,\tau)dx\,\nonumber\\
    &&\hspace{2cm}+\dfrac{\sigma_S}{2} m_I(t) \int_{\mathbb{R}_+}x^{2}\dfrac{d^2}{dx^2}\varphi_S(x)f_S(x,\tau)dx+\alpha m_R(t)\int_{\mathbb{R}_+} \dfrac{d}{dx}\varphi_S(x)f_S(x,\tau)dx\,.
\end{eqnarray}
Therefore, $f_S^\varepsilon(x,t)$ converges weakly to $f(x,t)$,  solution of the  Fokker-Planck equation
\begin{equation}
    \dfrac{\partial f_S(x,t) }{\partial t}=\dfrac{\sigma_S}{2}m_I(t) \dfrac{\partial^2}{\partial x^2}\left[x^{2}f_S(x,t)\right]+\dfrac{\partial}{\partial x}\left[\left(\beta m_I(t)x-\alpha m_R(t)\right)f_S(x,t)\right]\,,
    \label{FP_1}
\end{equation}
coupled with the no-flux boundary conditions
\begin{equation}
\label{eq:nfs}
\begin{split}
    &\left.\dfrac{\sigma_S}{2}m_I(t)\dfrac{\partial}{\partial x}\left[x^{2}f_S(x,t)\right]+\left[\beta m_I(t)x-\alpha m_R(t)\right]f_S(x,t)\right|_{x=0}=0\,,\\
    &\left. x^{2}f_S(x,t)\right|_{x=0}=0\,.
    \end{split}
\end{equation}
Analogous computations can be done for the kinetic equations \fer{weak_I} and \fer{weak_R}.  By doing so, in the quasi-invariant limit characterized by \fer{eq:scaling}, the system of Boltzmann type kinetic equations takes the form of a system of Fokker-Planck type equations with variable in time coefficients of diffusion and drift parameters. This system reads 
\begin{equation}\label{eq:system_PDE}
\begin{split}
\dfrac{\partial f_S(x,t)}{\partial t}&=\dfrac{\sigma_S}{2}m_I(t)\dfrac{\partial^2}{\partial x^2}\left[x^{2}f_S(x,t)\right]+\dfrac{\partial}{\partial x}\left[\left(\beta m_I(t)x-\alpha m_R(t)\right)f_S(x,t)\right]\,,\\
\dfrac{\partial f_I(x,t)}{\partial t}&=\dfrac{\sigma_I}{2} m_S(t)\dfrac{\partial^2}{\partial x^2}\left[x^{2}f_I(x,t)\right]+\dfrac{\partial}{\partial x}\left [(\left(\gamma\theta-\beta m_S(t)\right)x - \gamma(\theta-1)m_I) f_I(x,t)\right ]\,,\\
\dfrac{\partial f_R(x,t)}{\partial t}&= \dfrac{\sigma_R(\alpha)}{2}\dfrac{\partial^2}{\partial x^2}\left[x^2 f_R(x,t)\right] +  \dfrac{\partial}{\partial x} \left[(\alpha x - \gamma m_I(t))f_R(x,t) \right]\,.
\end{split}
\end{equation}
The first equation in system \fer{eq:system_PDE} is coupled with the no-flux boundary conditions defined in \eqref{eq:nfs}. Analogously, the second equation for the density $f_I$ is coupled with the no-flux boundary conditions
\begin{equation}
\label{eq:nfi}
\begin{split}
	&\left.\dfrac12\sigma_Im_S(t)\dfrac{\partial}{\partial x}\left[x^2f_I(x,t)\right]+\left[(\left( \gamma\theta - \beta m_S(t)\right)x - \gamma(\theta-1)m_I)\right]f_I(x,t)\right|_{x=0}=0\,,\\
	&\left. x^{2}f_I(x,t)\right|_{x=0}=0,
\end{split}
\end{equation}
Finally  the no-flux boundary conditions for the density $f_R$ read
\begin{equation}
\label{eq:nfr}
\begin{split}
	\dfrac{\sigma_R(\alpha)}{2}\dfrac{\partial}{\partial x}[x^2f_R(x,t)] + (\alpha x-\gamma m_I(t)) f_R(x,t) \Big|_{x = 0} = 0\,, \\
	x^2 f_R(x,t)\Big|_{x=0} = 0\,.
	\end{split}
\end{equation}
 \begin{remark} It is important to outline that the  Fokker--Planck type equations of system \fer{eq:system_PDE}, while characterized in general by variable in time coefficients of diffusion and drift, exhibit the same structure of a Fokker--Planck introduced in the framework of wealth distribution in a western society \cite{bouchaud2000,cordier2005}. In its standard form, this equation reads
 \be\label{eq:wealth}
 \dfrac{\partial f(x,t)}{\partial t}=\dfrac{\sigma}{2}\dfrac{\partial^2}{\partial x^2}\left[x^{2}f(x,t)\right]+ \dfrac{\partial}{\partial x}\left[(\lambda x- \mu)f(x,t)\right],
 \ee
 where the positive constants $\lambda$ and $\sigma$ are related to the details of the economic trades.
 In reason of this analogy, many results about existence, uniqueness and asymptotic behavior of the solutions to system \fer{eq:system_PDE} follow by a suitable generalization of the analogous ones proven for  equation \fer{eq:wealth} \cite{torregrossa2018}. As shown in \cite{torregrossa2018}, an important tool to work with equation \fer{eq:wealth}, and consequently with system  \fer{eq:system_PDE}, is the passage to Fourier transforms, where, as usual, for a given probability density $g(x)$ defined on $\R_+$, its Fourier transform is defined as
 \be\label{eq:Fourier}
 \hat g(\xi) = \int_{\R_+} g(x) e^{-i\xi x}\, dx.
 \ee
 For the sake of completeness, we will briefly detail some of these results in Appendix \ref{sec:wealth}.
 \end{remark}

\subsection{Moments of the simplified model}

At variance with the kinetic system \fer{eq:kine}, the Fokker--Planck type equations \eqref{eq:system_PDE} exhibit a number of interesting properties.  It is immediate to verify  that the evolution of the mean values coincides with the one defined in \eqref{eq:mSIR},  and corresponds to the classical SIR model. 

\begin{figure}
\centering
\includegraphics[scale = 0.2]{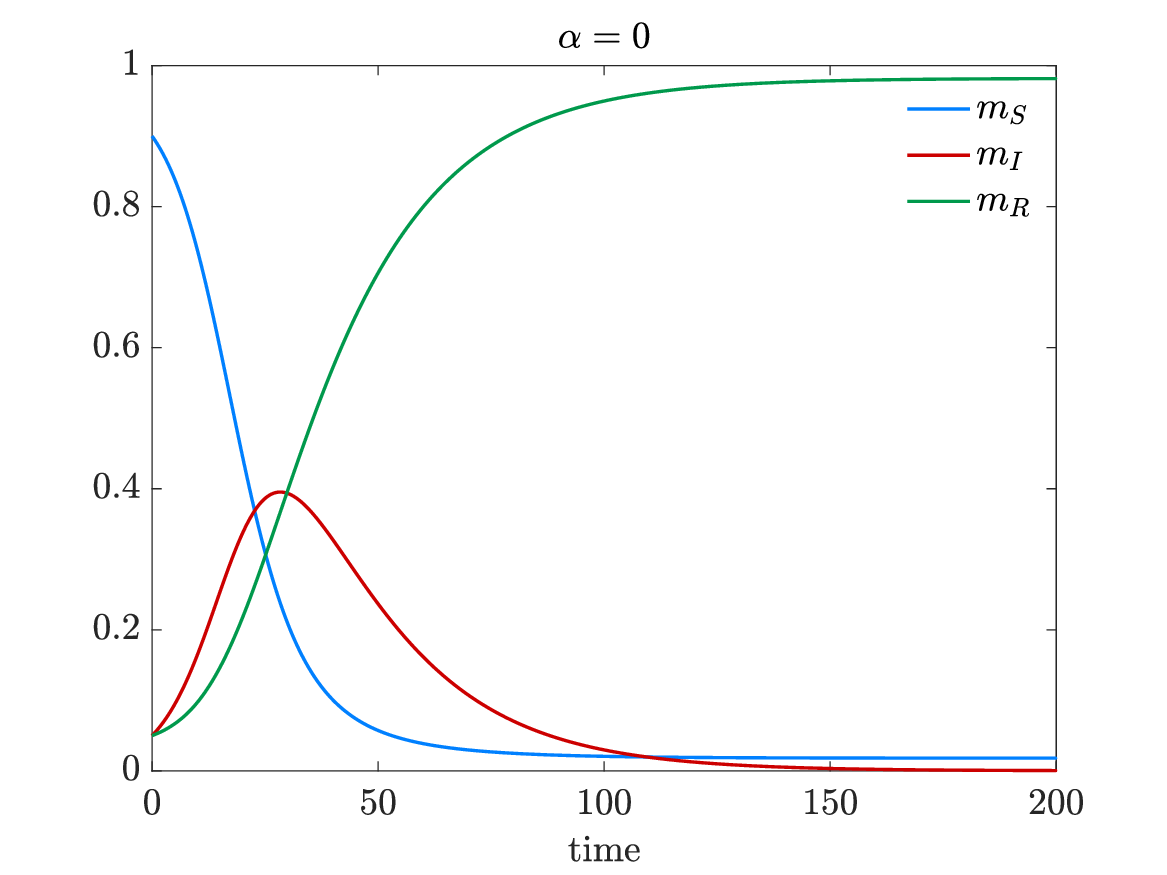}
\includegraphics[scale = 0.2]{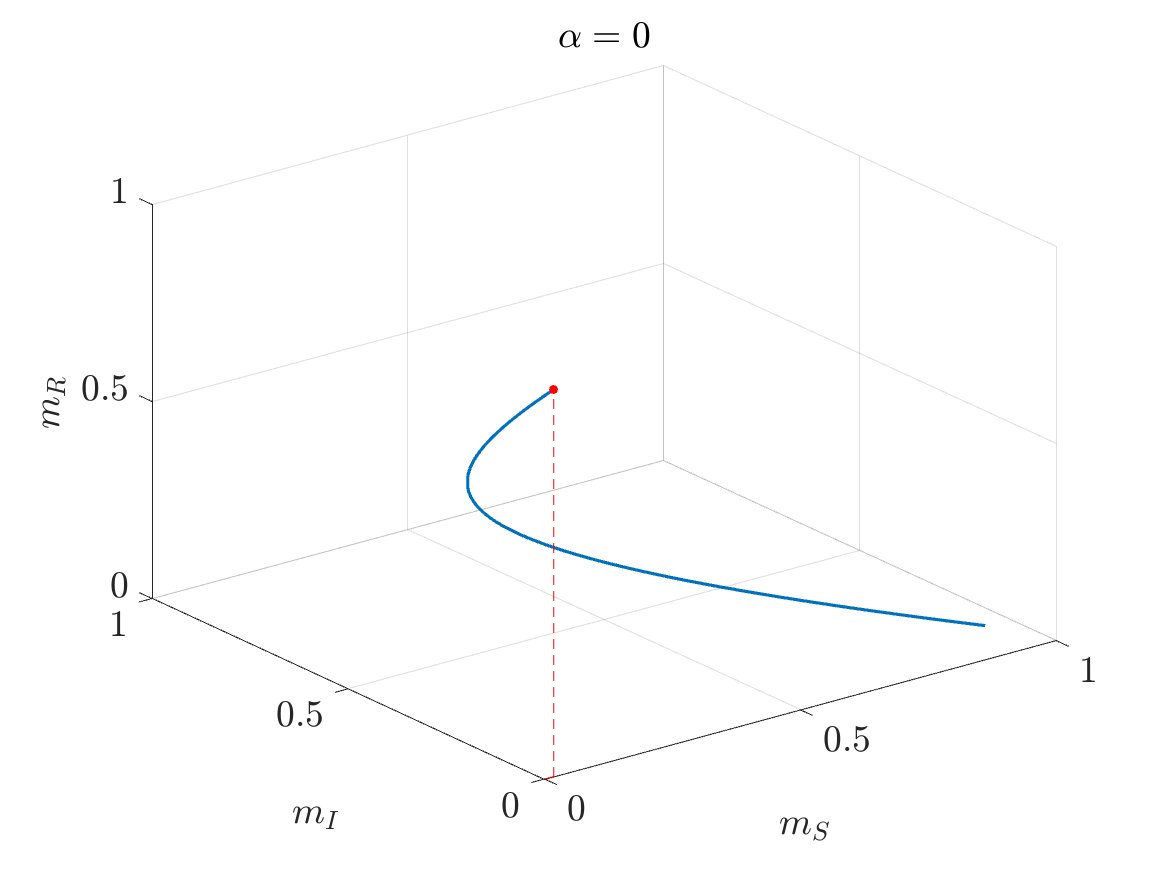}
\includegraphics[scale = 0.2]{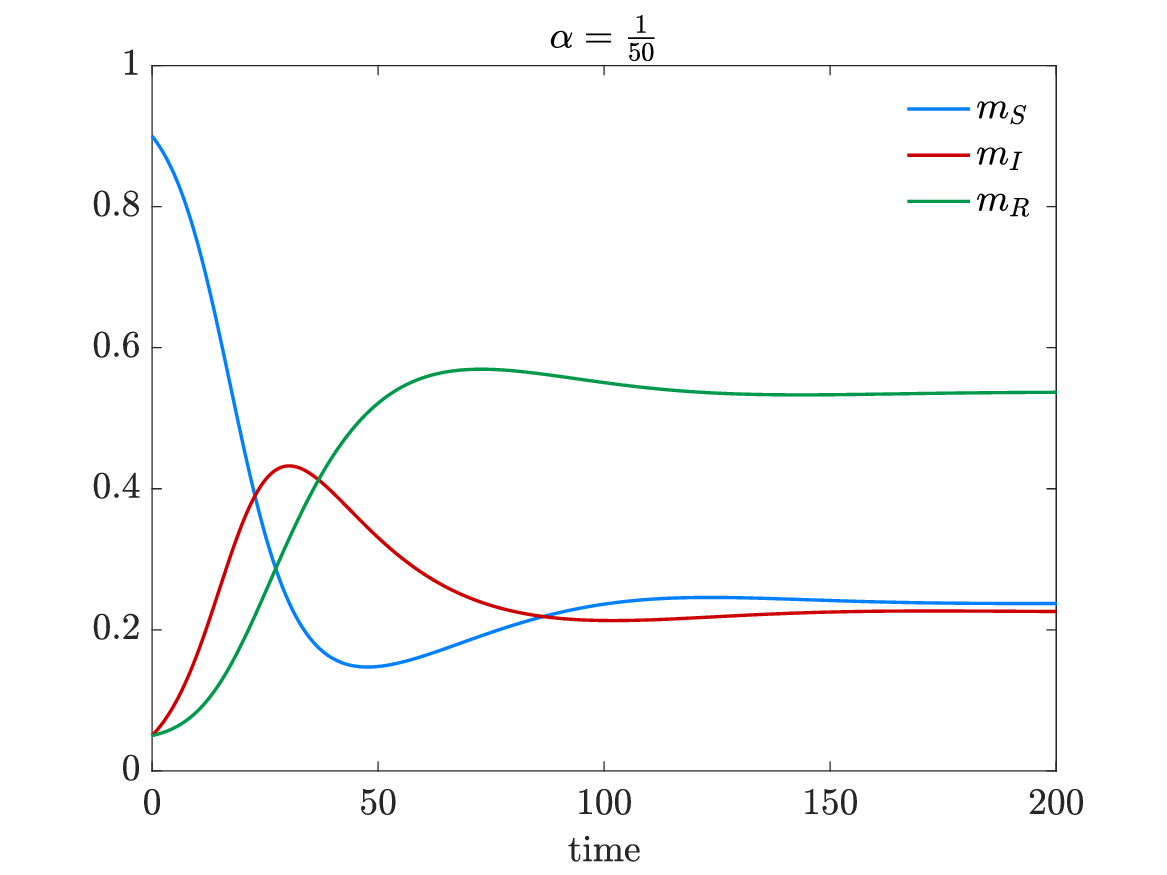}
\includegraphics[scale = 0.2]{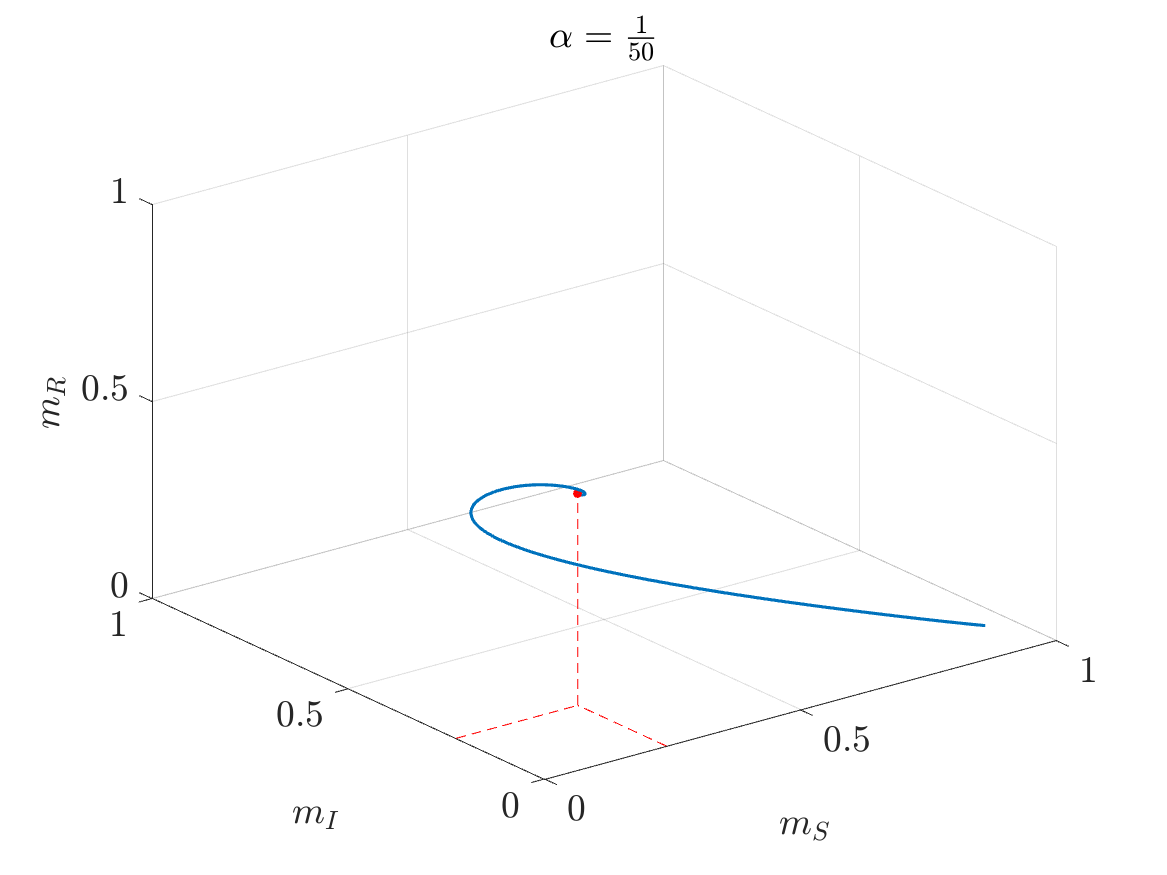}\\
\includegraphics[scale = 0.2]{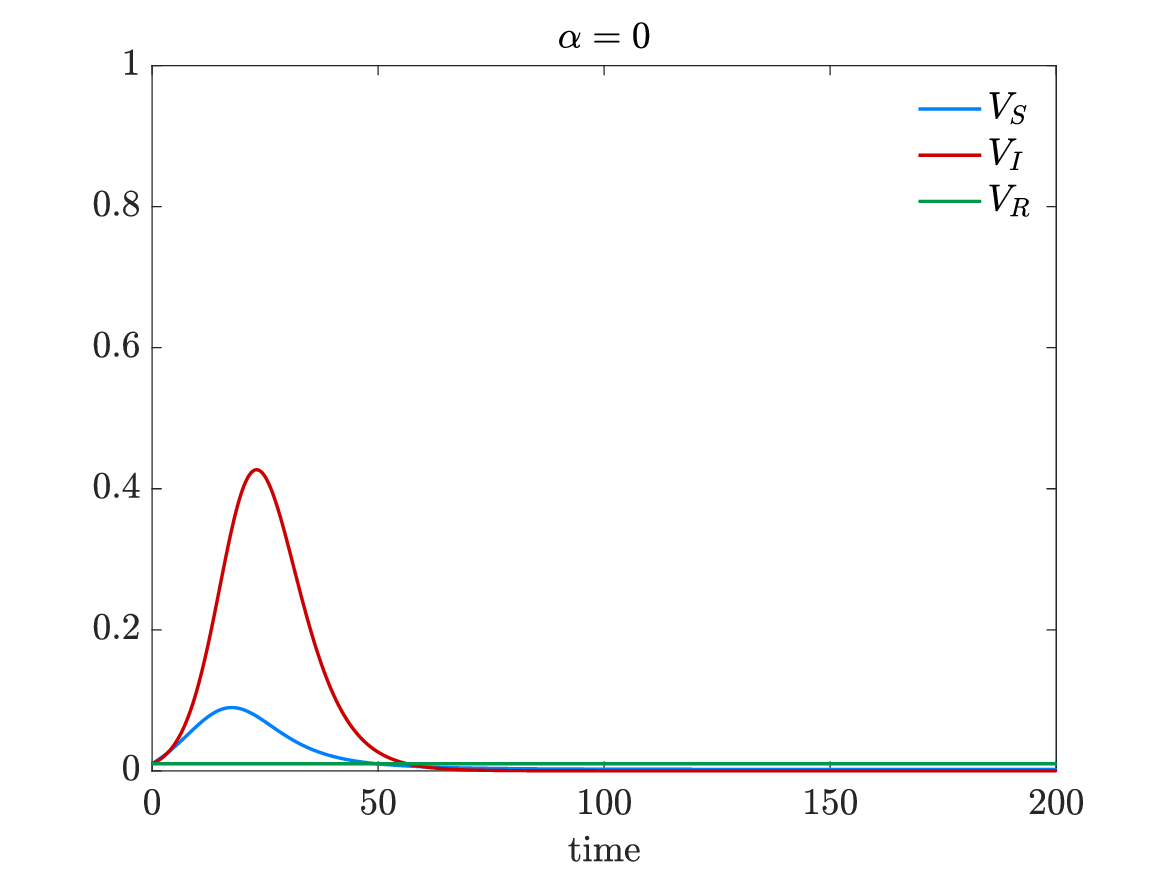}
\includegraphics[scale = 0.2]{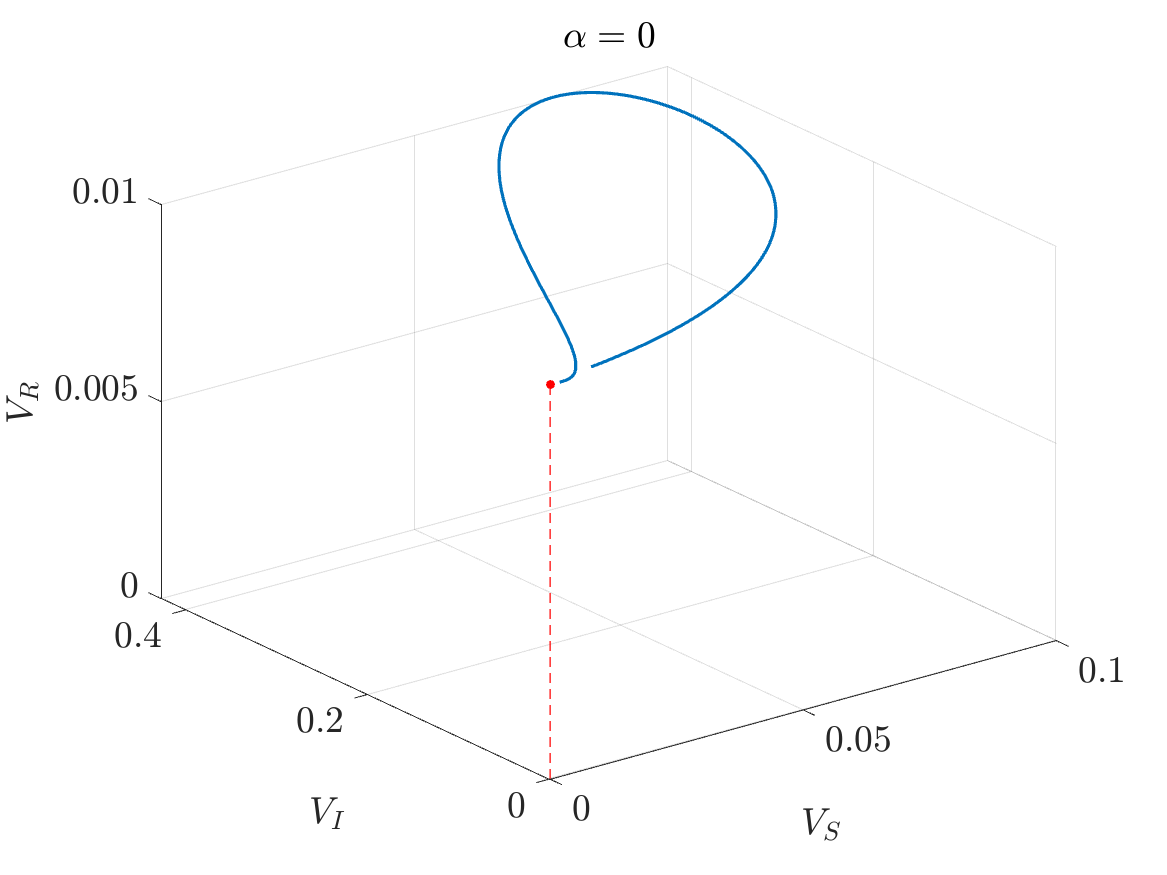}
\includegraphics[scale = 0.2]{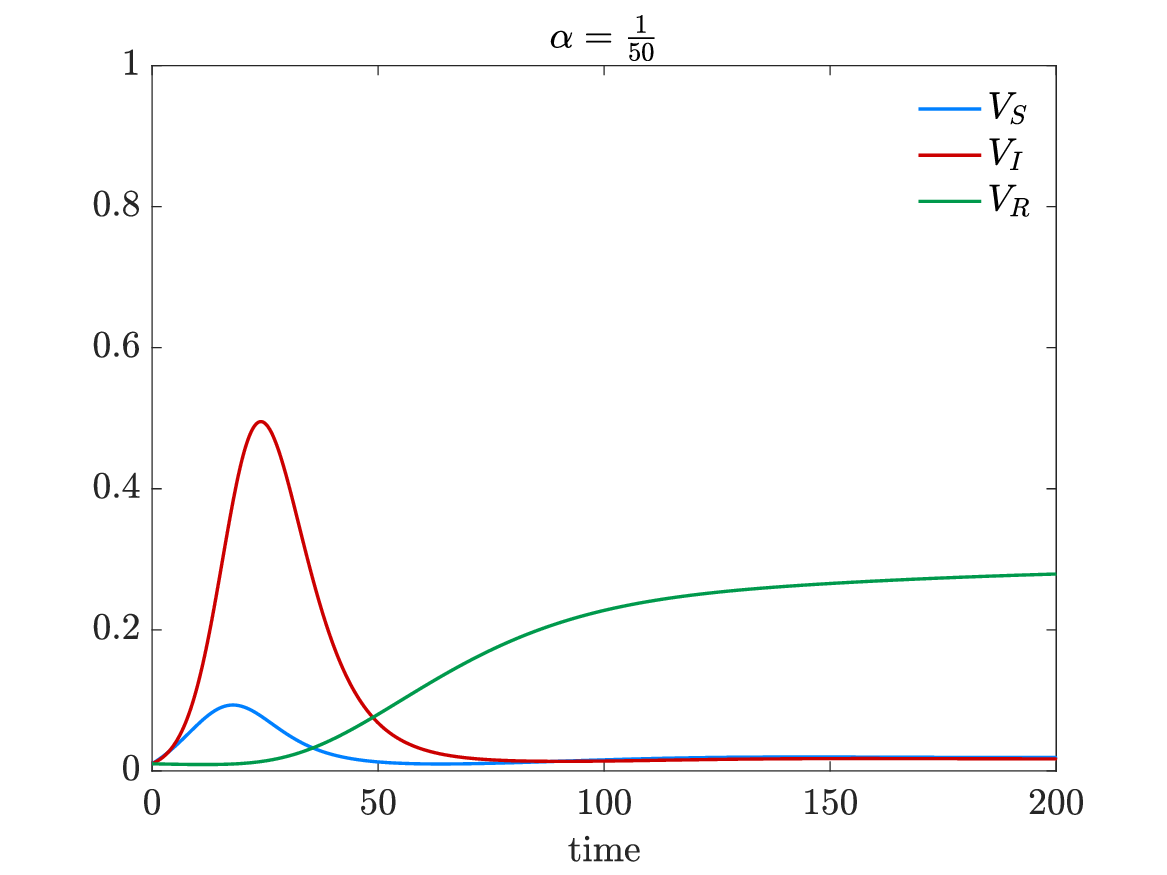}
\includegraphics[scale = 0.2]{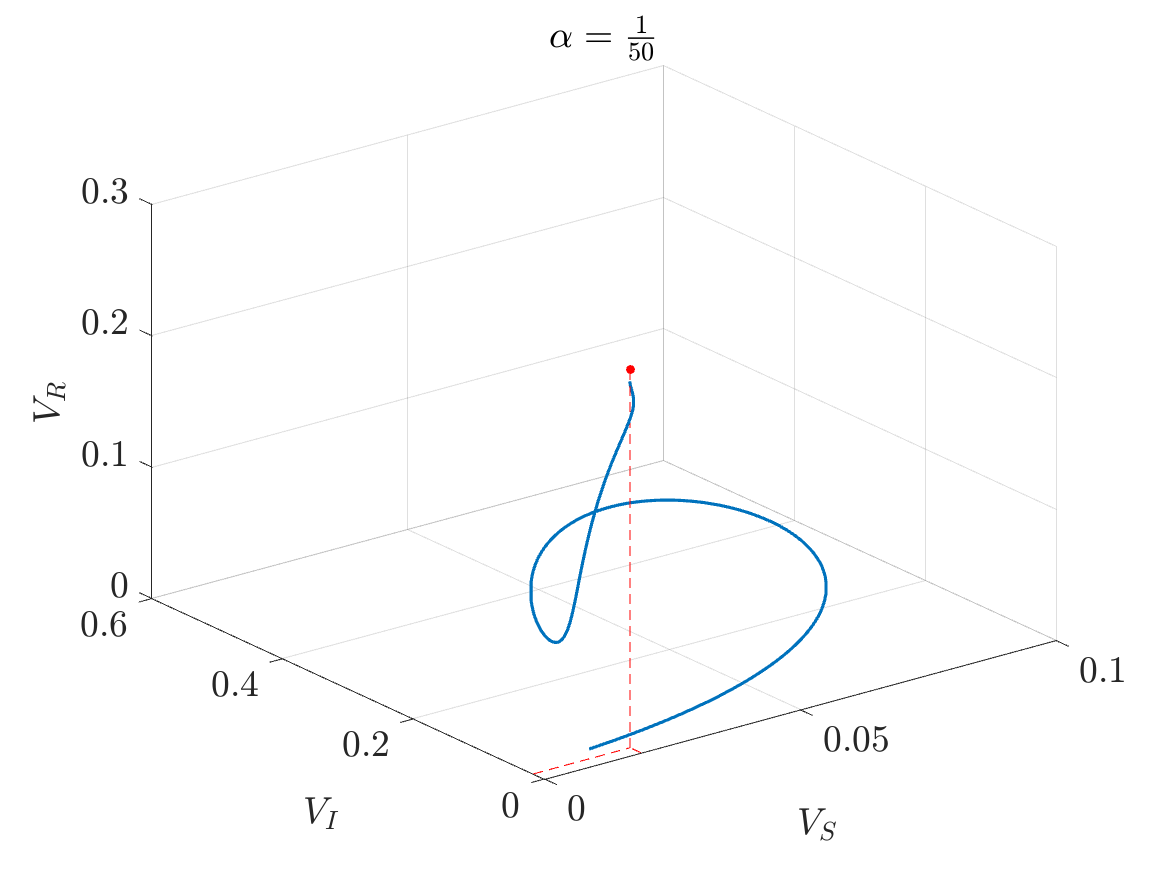}\\
\caption{Evolution of $m_J(t)$ and dynamics towards equilibrium (top row). Evolution of $V_J(t)$ and dynamics towards equilibrium point (bottom row). Parameters: $\beta = 2/10$, $\gamma = 1/21$, $\theta = 2$, $\sigma = 1/10$, $m_S(0) = 0.9$, $m_I(0) = m_R(0) = 0.05$, $v_J(0) = 0.01$, $J \in\mathcal{C}$.  }
\label{fig:0}
\end{figure}

Furthermore, owing to \fer{var-ev}, the evolution of the variances of the solutions to system \fer{eq:system_PDE} is given by 
\begin{equation}
\label{eq:variance_SIR_PDE}
\begin{split}
\dfrac{d}{dt}V_S(t) &= (\sigma_S -2\beta)m_I(t) V_S(t) + \sigma_Sm_I m_S^2(t)\,,  \\
\dfrac{d}{dt} V_I(t) &= \left[(\sigma_I+ 2\beta) m_S(t) -2\gamma\theta\right] V_I(t) +\sigma_Im_S(t) m_I^2(t) \,, \\
\dfrac{d}{dt}V_R(t) &= ( \sigma_R(\alpha) -2\alpha) V_R(t) + \sigma_R(\alpha) m_R(t)^2\,.
\end{split}
\end{equation}
In this case, provided  {$\sigma_S < 2\beta$,  {$\sigma_I < 2\beta(\theta-1) $}  and $\sigma_R(\alpha) < 2\alpha$}, the variances remain uniformly bounded in time, and converge towards explicit values,  given by
\begin{equation}\label{eq:Vinfty}
(V_S^{\infty}, V_I^{\infty}, V_R^{\infty} ) =
\begin{cases}
\displaystyle  \left(\frac{\sigma_S (m_S^\infty)^2}{2\beta - \sigma_S}, \frac{\sigma_I m_S^\infty (m_I^\infty)^2}{2\gamma\theta-(\sigma_I +2\beta)m_S^\infty},{\frac{\sigma_R(\alpha)}{2\alpha- \sigma_R(\alpha)}(m_R^\infty)^2}\right) & \alpha>0 \\
 \displaystyle \left(\frac{\sigma_S( m_S^\infty)^2}{2\beta - \sigma_S},0,V_R(0)\right) & \alpha = 0
\end{cases}\,.
\end{equation}
We can notice how, if $\alpha >0$, the three variances of the  populations remain bounded away from zero.  This suggests that the solutions to system \fer{eq:system_PDE} should converge towards an asymptotic profile. At variance, when $\alpha =0$, while the variance of the susceptible and recovered populations converge to  finite values,  the variance of  the population of  infected individuals annihilates asymptotically. Thus, in the case without re-infection, the epidemic spreading will end up in absence of infected individuals. In Figure \ref{fig:0} we depict the evolution of the system \eqref{eq:mSIR}-\eqref{eq:variance_SIR_PDE} together with the equilibrium points. 

\begin{remark}\label{rem:var}
It is important to notice that, under the conditions {$\sigma_S < 2\beta$ and $\sigma_R(\alpha) < 2\alpha$}, the coefficients of $V_S(t)$ and $V_R(t)$ in  \fer{eq:variance_SIR_PDE} are negative, while, under the condition $\sigma_I < 2\beta(\theta-1) $, the coefficient of $V_I(t) $ in the second  equation  \fer{eq:variance_SIR_PDE} is asymptotically negative, due to the asymptotic value $m_S^\infty = \gamma/\beta$. This clearly implies that there exists a time $\bar t$ such that, for $t \ge \bar t$, $(\sigma_I+ 2\beta) m_S(t) -2\gamma\theta <0$, and $V_I(t)$ stays uniformly bounded in time.
\end{remark}

\begin{remark}
The  evolution of the variances in \eqref{eq:variance_SIR_PDE} can be obtained from the evolution of the variances of the Boltzmann-type model \eqref{eq:VS_B}-\eqref{eq:VI_B}-\eqref{eq:VR_B} under the quasi-invariant scaling \eqref{eq:scaling}. 
\end{remark}

\begin{remark}
We can observe how, if  $\alpha = 0$, the asymptotic mean number of infected converges to zero triggering the convergence of $V_I^{\infty}$ to zero. In this case we may expect to obtain a Dirac delta distribution for the infected population. 
\end{remark}

\subsection{Local equilibria}\label{sec:local}
Detailed informations about the distributions of the populations in time can be recovered by looking at the so-called \emph{quasi-equilibrium states} \cite{BonMenTosZan}. 
The quasi-equilibrium states of the Fokker–Planck system \fer{eq:system_PDE} are the density
functions $f_J^{q}(x, t)$, $J \in\mathcal{C}$ solutions to the  differential equations
\begin{equations}\label{eq:q-e}
    &\dfrac12\sigma_S\dfrac{\partial}{\partial x}\left[x^{2}f^q_S(x,t)\right]+\left(\beta x-\alpha \frac{m_R(t)}{m_I(t)}\right)f^q_S(x,t)=0, \\
    &\dfrac12\sigma_I\dfrac{\partial}{\partial x}\left[x^{2}f^q_I(x,t)\right]+\left[\left(\frac{\gamma\theta -\beta m_S(t)}{m_S(t)}\right)x-\gamma(\theta-1)\frac{m_I(t)}{m_S(t)}\right]f^q_I(x,t)=0,  \\ 
  &\dfrac{\sigma_R(\alpha)}{2}\dfrac{\partial}{\partial x}\left[x^2 f^q_R(x,t)\right] + [(\alpha x - \gamma m_I(t))f^q_R(x,t)] =0\,.
\end{equations}
The solutions of \fer{eq:q-e} can be easily obtained resorting to expression \fer{equi2} of the Appendix, with the difference that now the coefficients $\sigma$, $\lambda$ and $\mu$ may depend on time.  In this case  the local equilibria $f_S^q(x,t)$ and $f_I^q(x,t)$ are inverse Gamma densities with time dependent shape parameters $\nu_S$ (respectively $\nu_I(t)$) and scale parameters $\omega_S(t)$ (respectively $\omega_I(t)$). The explicit form of the local equilibria is    \begin{equation}\label{local}
    \begin{split}
       & f_S^{q}(x,t)=\dfrac{\omega_S^{\nu_S}(t) }{\Gamma(\nu_S)}\left(\dfrac{1}{x}\right)^{\nu_S+1}e^{-\omega_S(t)/x}, \quad f_I^{q}(x,t)=\dfrac{\omega_I(t)^{\nu_I(t)}}{\Gamma(\nu_I(t))}\left(\dfrac{1}{x}\right)^{\nu_I(t)+1}e^{-\omega_I(t)/x}\,, \\
        & \nu_S=1 +\dfrac{2\beta}{\sigma_S}\,,\quad\omega_S(t)=\dfrac{2\alpha}{\sigma_S}\,\dfrac{m_R(t)}{m_I(t)}\,,\\
        &  \nu_I(t)=1+ \dfrac{2( \gamma\theta-\beta m_S(t))}{\sigma_I m_S(t)}\,,\quad\omega_I(t)=\dfrac{2\gamma(\theta- 1)m_I(t)}{\sigma_Im_S(t)}\,.
        \end{split}
    \end{equation}
    Furthermore, if $\alpha>0$ the solution to the third equation is the inverse Gamma density 
    \begin{equation}
    \label{eq:localR}
    f_R^q(x,t) = \dfrac{\omega_R^{\nu_R}}{\Gamma(\nu_R)}\left( \dfrac{1}{x}\right)^{\nu_R+1}e^{-\omega_R(t)/x},
    \end{equation}
    where the shape parameter $\nu_R = 1+2\alpha/\sigma_R(\alpha)$ and the scale parameter $\omega_R = \frac{2\gamma m_I(t)}{\sigma_R(\alpha)}>0$. 
 \begin{remark}
 It is important to remark that, since the local equilibria \fer{local} are inverse gamma densities,  the existence of moments is closely related to the values of the shape  parameters. In particular, the local equilibria have bounded mean values if and only if $\nu_S>1$ and $\nu_I(t) >1$. While this condition is always satisfied for $\nu_S$, for $\nu_I(t)$ it is required that
 \be\label{nu_I}
 \beta m_S(t) < \gamma\theta. 
 \ee
Since, thanks to \fer{eq:minfty}, this condition is asymptotically satisfied since $\theta>1$, we can easily conclude that the mean values of the local equilibria are uniformly bounded from above. Likewise, as shown in Remark \ref{rem:var},  the uniform boundedness of the variances, which requires $\nu_S>2$ and $\nu_I(t) > 2$,  is satisfied provided
\be\label{nu-var}
\sigma_S < 2\beta, \quad \sigma_I < 2\beta(\theta-1) .
\ee
In both cases, these conditions can be achieved in presence of small values of the diffusion constants $\sigma_S$ and $\sigma_I$, or in presence of a sufficiently big value of the constant $\theta$.
It is equally important to outline that, as observed in Remark \ref{rem:mean}, the uniform boundedness from below and above of the mean values implies the uniform boundedness from below and above  of the shape and scale parameters of the local equilibria $f_S^q(x,t)$ and $f_I^q(x,t)$, as well as their convergence to explicit limit values. {Last, since when $\alpha >0$ it holds $\sigma_R(\alpha) < 2\alpha$, $\nu_R >2$, and the variance of the distribution of the population of recovered is always bounded, independently of the value of $\alpha >0$}.
 \end{remark} 
 
 \section{Large-time behavior of the Fokker--Planck system}\label{sec:large}
 
 As discussed in Section \ref{sec:FP}, the convergence of the means of the solutions to system \fer{eq:system_PDE} towards explicit values, coupled with the fact that at any time $t >0$ the local equilibria, in the form of inverse Gamma densities, have shape and scale parameters that are explicitly expressed in terms of these time-dependent mean values, so that they converge towards inverse Gamma densities with constant parameters,   suggests that the solutions to  system \fer{eq:system_PDE} should converge towards the asymptotic values of the local equilibria, both in the case $\alpha >0$ and in the case $\alpha = 0$.  In this Section, we will show that this behavior can be fully justified.  However the cases $\alpha >0$  and $\alpha =0$ lead to completely different results. While in presence of re-infection, corresponding to $\alpha >0$, the asymptotic profiles lose memory of the initial distributions of the populations, in absence of re-infection, so that $\alpha =0$, the asymptotic profiles retain memory of the initial configuration. 
 
 \subsection{The large-time behavior in the case with re-infection}
 
In the rest of the paper, we will assume that the conditions on parameters specified in Remark  \ref{rem:var}
\begin{equation*}
	\sigma_S < 2\beta\,,\qquad\sigma_R(\alpha) < 2\alpha\,,\qquad\sigma_I < 2\beta(\theta-1)
\end{equation*}
are satisfied.  
In this case, it can be easily verified that, since the time derivatives of local equilibria are explicitly related to the local equilibria themselves, the $L_1(\R_+)$-norm of these time derivatives can be expressed as a linear combination of the time derivatives of the shape and scale parameters with uniformly bounded coefficients. 

Indeed, for a given inverse Gamma distribution $h(x,t)$ of shape parameter $\nu(t)>2$ and scale parameter $\omega(t)$, it holds
\begin{equations}\label{deri}
&\frac{\partial \log h(x,t)}{\partial t} = \frac{\partial }{\partial t} \left[ \nu(t)\log\omega(t) - \log \Gamma(\nu(t))- [1+\nu(t)]\log x -\frac{\omega(t)}x \right] = \\
&\nu^\prime(t)\left[ \log\omega(t) -\psi(\nu(t)) -\log x\right] + \omega^\prime(t)\left[  \frac{\nu(t)}{\omega(t)}  - \frac1x\right],
\end{equations}
where $\psi(z) = \frac d{dz}\log \Gamma(z)$ is the $\psi$-function \cite{GR}, that for $z >0$ reads
\[
\psi(z) = \log z -\frac 1{2z} - 2 \int_0^{+\infty} \frac{\tau \, d\tau}{(\tau^2+ z^2)(e^{2\pi\tau}-1)}.
\]
Consequently
\be\label{bene}
\frac{\partial  h(x,t)}{\partial t} =\left\{ \nu^\prime(t)\left[ \log\omega(t) -\psi(\nu(t)) -\log x\right] + \omega^\prime(t)\left[  \frac{\nu(t)}{\omega(t)}  - \frac1x\right]\right\}h(x,t).
\ee
Now, consider that 
\begin{equations}\nonumber
&|\log x| = |x\log x|\frac 1x \le \frac 1{ex}, \quad 0< x \le 1\,, \\
& \log x \le 1+ x , \qquad x >1.
\end{equations}
Therefore 
\be\label{bound-h}
\left|\frac{\partial  h(x,t)}{\partial t}\right| \le \left\{ |\nu^\prime(t)|\left[ |\log\omega(t) -\psi(\nu(t))|  + 1 +x + \frac 1{ex}\right] + |\omega^\prime(t)|\left[  {\frac{\nu(t)}{\omega(t)}}  + \frac1x\right]\right\}h(x,t).
\ee
On the other hand, for the inverse gamma distribution $h(x,t)$ it holds
\be\label{means} 
\int_{\R_+} x h(x,t) \, dx = \frac{\omega(t)}{\nu(t) -1}\,, \quad \int_{\R_+} \frac 1x h(x,t) \, dx = \frac{\nu(t)}{\omega(t) }\,.
\ee
The upper bound \fer{bound-h}, coupled with \fer{means} implies that 
\be\label{stima}
\int_{\R_+} \left|\frac{\partial  h(x,t)}{\partial t}\right|\, dx \le |\nu^\prime(t)| A(t) +  |\omega^\prime(t)|B(t),
\ee
where
\begin{equations}\label{coeff1}
A(t) &= |\log\omega(t) -\psi(\nu(t))|  + 1 + \frac{\nu(t)}{e \omega(t)} +\frac{\omega(t)}{\nu(t) -1}\,, \qquad
B(t) = 2\,\frac{\nu(t)}{\omega(t)}\,.
\end{equations}
If now we assume that  $c_1 \le \nu(t) \le C_1$ and $c_2 \le \omega(t) \le C_2$, it follows that $A(t) \le D_A$ and $B(t) \le D_B$, where $D_A$ and $D_B$ are suitable bounded constants.

Making use of the previous results,  we can now proceed to study the large-time behavior of the densities $f_J(x,t)$, $J\in \mathcal{C}$. 

We start our analysis by considering first the case of the density $f_S(x,t)$ of susceptible individuals, that satisfies the first Fokker--Planck equation in \fer{eq:system_PDE}.

 For any given time $t >0$, the local equilibrium of this Fokker--Planck equation with variable in time coefficient of diffusion and drift is the (unique) function  $f_S^q(x,t)$  of unit mass satisfying 
\begin{equation}\label{eq:local}
\dfrac{\sigma_S}{2}m_I(t)\dfrac{\partial}{\partial x}\left[x^{2}f_S^q(x,t)\right]+\left(\beta m_I(t)x-\alpha m_R(t)\right)f_S^q(x,t) = 0\,.
\end{equation}
Therefore, since the right-hand sides of the (linear) Fokker--Planck operators in  \fer{eq:system_PDE} vanish when evaluated in correspondence to the quasi-equilibrium densities, we have the identity
\begin{equations}
&\dfrac{\partial [f_S(x,t)- f_S^q(x,t)]}{\partial t}= - \dfrac{\partial f_S^q(x,t)}{\partial t} +\\
& \dfrac{\sigma_S}{2}m_I(t)\dfrac{\partial^2}{\partial x^2}\left\{x^{2}[f_S(x,t)- f_S^q(x,t) \right\}+\dfrac{\partial}{\partial x}\left\{\left(\beta m_I(t)x-\alpha m_R(t)\right)[f_S(x,t)-f_S^q(x,t)]\right\}\,.
\end{equations}
Let $1/2 <p < 1$. By passing to Fourier transform, and applying the  procedure described in Appendix \ref{sec:wealth}, we easily obtain that the $\dot H_{-p}$--Sobolev norm of the difference $\fc_S(\xi,t)- \fc_S^q(\xi,t)$ satisfies the inequality 
\begin{equations}\label{inizio}
&\dfrac d{dt}{\int_\R |\xi|^{-2p} |\fc_S(\xi,t)- \fc_S^q(\xi,t)|^2}\,d\xi \le  \int_\R  |\xi|^{-2p}\left|\widehat{ \dfrac{\partial f_S^q}{\partial t}}(\xi,t)\right||\fc_S(\xi,t)- \fc_S^q(\xi,t)|\, d\xi +\\
& -(2p-1)\left[m_S(t)\left(\sigma_I \frac{3-2p}4 + \beta \right)  \right] {\int_\R |\xi|^{-2p} |\fc_S(\xi,t)- \fc_S^q(\xi,t)|^2}\,d\xi.
\end{equations}
Thanks to the bound \fer{stima} we have
\begin{equation}\label{eq:fq}
\left|\widehat{ \dfrac{\partial f_S^q}{\partial t}}(\xi,t)\right| =\left|  \int_{\R_+} \frac{\partial  f_S^q(x,t)}{\partial t}e ^{-i\xi x} \, dx \right| \le  
 \int_{\R_+} \left|\frac{\partial  f_S^q(x,t)}{\partial t}\right| \, dx \le \epsilon_S(t) ,
\end{equation}
where we set
\be\label{con11}
\epsilon_S(t) = |\nu_S^\prime(t)| A_S +  |\omega_S^\prime(t)|B_S. 
\ee
In \fer{con11}  $A_S$ and $B_S$ are the (finite) upper bounds of $A_S(t)$, and, respectively, of $B_S(t)$. 

Moreover, since both $f_S(x,t)$ and $f_S^q(x,t)$ are probability densities, so that $\fc_S(0,t) = \fc_S^q(0,t) =1$, expanding in Taylor series around the point $\xi =0$ up to the first order, we get
\be\label{eq:diff}
|\fc_S(\xi,t)- \fc_S^q(\xi,t)| = \left|\widehat{ \dfrac{\partial f}{\partial \xi}}(\xi^*,t)- \widehat{ \dfrac{\partial f_S^q}{\partial \xi}}(\bar \xi,t) \right| |\xi| \le (m_S(t) + m_S^q(t))|\xi| \le M_S|\xi|.
\ee
Indeed, for a given probability density $h(x)$ of finite mean value, integration by parts gives
\[
\left|\widehat{ \dfrac{\partial h(\xi)}{\partial \xi}}\right| = \left|  \int_{\R_+} \frac{\partial  h(x)}{\partial x}e ^{-i\xi x} \, dx\right| = \left| i \int_{\R_+} x\, h(x)e ^{-i\xi x} \, dx\right| \le  \int_{\R_+} x\, h(x)\, dx.
\]
Therefore, for any given $\xi$
\be\label{Rmeno}
 |\xi|^{-2p}\left|\widehat{ \dfrac{\partial f_S^q}{\partial t}}(\xi,t)\right||\fc_S(\xi,t)- \fc_S^q(\xi,t)| \le  M_S\epsilon_S(t) |\xi|^{1-2p}\,.
\ee
Making use of \fer{Rmeno}, for any given $R >0$ we obtain
\be\label{int-meno}
\int_{|\xi| \le R}  |\xi|^{-2p}\left|\widehat{ \dfrac{\partial f_S^q}{\partial t}}(\xi,t)\right||\fc_S(\xi,t)- \fc_S^q(\xi,t)| \, d\xi \le M_S \epsilon_S(t) \frac  {R^{2-2p}}{1-p}. 
\ee
On the other hand, by Cauchy-Schwarz inequality 
\begin{equations}\label{Rpiu}
& \int_{|\xi| > R}  |\xi|^{-2p}\left|\widehat{ \dfrac{\partial f_S^q}{\partial t}}(\xi,t)\right||\fc_S(\xi,t)- \fc_S^q(\xi,t)| \, d\xi \le \\
&\left(   \int_{|\xi| > R} |\xi|^{-2p}\left|\widehat{ \dfrac{\partial f_S^q}{\partial t}}(\xi,t)\right|^2 \, d\xi\right)^{1/2}\left( \int_{|\xi| > R}  |\xi|^{-2p}|\fc_S(\xi,t)- \fc_S^q(\xi,t)|^2 \, d\xi \right)^{1/2} \le \\
& \epsilon_S(t) \left(   \int_{|\xi| > R} |\xi|^{-2p}\, d\xi\right)^{1/2}\left( \int_{\R}  |\xi|^{-2p}|\fc_S(\xi,t)- \fc_S^q(\xi,t)|^2 \, d\xi \right)^{1/2} =\\
& \epsilon_S(t) \frac 1 { R^{p-1/2}}\left(\frac 2{2p-1} \int_{\R}  |\xi|^{-2p}|\fc_S(\xi,t)- \fc_S^q(\xi,t)|^2 \, d\xi \right)^{1/2}\,.
\end{equations}
Optimizing over $R$ we get 
\be\label{finale}
\int_{\R}  |\xi|^{-2p}\left|\widehat{ \dfrac{\partial f_S^q}{\partial t}}(\xi,t)\right||\fc_S(\xi,t)- \fc_S^q(\xi,t)| \, d\xi \le \epsilon_S(t) C_p M_S^{\frac{1}{3-2p}}\left( \int_{\R}  |\xi|^{-2p}|\fc_S(\xi,t)- \fc_S^q(\xi,t)|^2 \, d\xi \right)^{\frac{2-2p}{3-2p}},
\ee
where $C_p$ is an explicitly computable constant.

Let us set
\[
z(t) = \int_{\R}  |\xi|^{-2p}|\fc_S(\xi,t)- \fc_S^q(\xi,t)|^2 \, d\xi. 
\]
Then, thanks to \fer{inizio} and to \fer{finale} we have that $z(t)$ satisfies the differential inequality
\[
\frac{dz(t)}{dt} \le -c_pm_I(t)  z(t) + \epsilon_S(t) C_p M_S^{\frac{1}{3-2p}}z(t)^{\frac{2-2p}{3-2p}},
\]
where the value of the constant $c_p$ can be recovered from \fer{stima}. Since $z(t) >0$, by setting $y(t) = z(t)^{1/(3-2p)}$ it is immediate to show that this inequality is equivalent to  
\be\label{diffe}
\frac{dy(t)}{dt} \le -\frac{c_pm_I(t)}{3-2p}  y(t) + \frac{C_p\epsilon_S(t) M_S^{\frac{1}{3-2p}}}{3-2p},
\ee
whose solution, since $\epsilon_S(t) \to 0$ as $t \to + \infty$ converges to zero. 

Indeed, if we set $a(t) = c_pm_I(t)/(3-2p)$, and $b(t)= C_p\epsilon_S(t) M_S^{\frac{1}{3-2p}}/(3-2p)$, the explicit solution to equation \fer{diffe} reads
\[
y(t) = \left[ y(0) + \int_0^t b(s)  \exp\left\{\int_0^s a(\tau)\, d\tau \right\} \, ds \right] \exp\left\{-\int_0^t a(s)\, ds \right\}\,.
\]
Since $a(t) \ge c>0$, so that $\lim_{t\to\infty} \int_0^t a(s) \, ds = +\infty$, and $b(t) \to 0$, the limit value of $y(t)$ can be obtained by applying the l'Hopital rule. One obtains
\[
\lim_{t\to \infty} y(t) = \lim_{t\to \infty} \frac{b(t)}{a(t)} =0.
\]
The same analysis applies to the densities $f_I(x,t)$ of infectious individuals and $f_R(x,t)$ of recovered individuals.

The previous results allow to  determine both the shape of the three equilibrium  densities, and  the role played by the local equilibria, which, at least for intermediate times, provide a good approximation of the true solution of system \fer{eq:system_PDE}.  It is important to remark that, in presence of re-infection, the final state of the system completely loses memory of the initial distribution of the populations of susceptible, infectious and recovered individuals. On the contrary, as we shall see in the next Section \ref{sec:no-re-infection}, in absence of re-infection the population of susceptible individuals retains memory of the initial state. Once again, we shall confirm and extend known results; in particular, while it can be observed that the memory of the initial data is retained through equation \eqref{constant_motion}, which provides a relation between the initial and the asymptotic states, we shall also show that this dependence emerges in the distribution function and is thus reflected in each of its moments.

\subsection{Large--time behavior in absence of re-infection}\label{sec:no-re-infection}

A somewhat different behavior can be observed in absence of re-infection, corresponding to the choice $\alpha=0$. In this case,  the evolution of the mean values read
\begin{equations}\label{macro_1_red}
\frac{dm_S(t)}{dt}&= -\beta{m_I(t) m_S(t)}, \qquad m_S^0 = m_S(t=0) \ge 0,\\
\frac{dm_I(t)}{dt}&= \beta{m_I(t) m_S(t)} - \gamma m_I(t), \qquad m_I^0= m_I(t=0)  \ge 0, \\
\frac{dm_R(t)}{dt}&=  \gamma m_I(t), \qquad m_R^0 = m_R(t=0) \ge 0\,.
\end{equations}
Since $\alpha=0$, the Fokker--Planck equation for the susceptible's individuals in \fer{eq:system_PDE} takes the (simple) form
\begin{equation}
	\dfrac{\partial f_S(x,t)}{\partial t}= m_I(t) \left[\dfrac12\sigma_S\dfrac{\partial^2}{\partial x^2}\left[x^{2}f_S(x,t)\right]+\dfrac{\partial}{\partial x}\left(\beta xf_S(x,t)\right) \right]\,.
	\label{FP_1_red}
\end{equation}
A further simplification is obtained by setting  $y(x,t)=m_S(t)x$, and
\begin{equation}
	g_S(x,t)=m_S(t)f_S(y(x,t),t).
\end{equation}
It is immediate to realize that the mean value of the scaled density $g_S(x,t)$ remains  equal to $1$ in time. Furthermore, the density $g_S(x,t)$ satisfies a pure diffusion equation. This equation is well-known to people working in probability theory and finance, since it describes a particular version of the geometric Brownian motion \cite{tos2016}.  This result can be easily verified by direct computations. Indeed
\begin{equations}\label{semply}
\dfrac{\partial}{\partial t}\,g_S(x,t)&= \dot{m}_S(t)f_S(y,t)+m_S(t)\dfrac{\partial}{\partial y}f_S(y,t)\dfrac{\partial y}{\partial t}+m_S(t)\dfrac{\partial}{\partial t}f_S(y,t) \\
	&=\dfrac{\dot{m}_S(t)}{m_S(t)}\left[g_S(x,t)+x\dfrac{\partial}{\partial x}g_S(x,t)\right]+m_S(t)\dfrac{\partial}{\partial t}f_S(y,t)  \\
	&=-\beta m_I(t) \dfrac{\partial}{\partial x}[xg_S(x,t)]+m_S(t)\dfrac{\partial}{\partial t}f_S(y,t) ,
\end{equations}
where in last line we used the first equation in \fer{macro_1_red}. 
Note however that the last term can be rewritten by resorting to \eqref{FP_1_red}. We obtain
\begin{equations}
	\dfrac{\partial}{\partial t}\,g_S(x,t)&= -\beta m_I(t) \dfrac{\partial}{\partial x}\left(xg_S(x,t)\right)+\dfrac12\sigma_Sm_I\dfrac{\partial^2}{\partial x^2}\left(x^2g_S(x,t) \right)+\beta m_I(t)\dfrac{\partial}{\partial x}\left(xg_S(x,t)\right)\\
	&=\dfrac12\sigma_Sm_I(t)\dfrac{\partial^2}{\partial x^2}\left(x^2g_S\right).
\end{equations}
 Consequently, the density $g(x,t)$ satisfies a geometric Brownian motion with time-dependent coefficient 
\[
\dfrac{\partial}{\partial t} g_S(x,t) = m_I(t)\dfrac{\sigma_S }{2} \dfrac{\partial^2}{\partial x^2} (x^2g_S(x,t)).
\]
We can easily reduce this equation to a classical geometric Brownian motion by setting $g_S(x,t) = u(x, t^*)$, where
\[
t^* = t^*(t) = \frac 12\sigma_S\int_0^t m_I(s) \, ds. 
\]
Indeed, since $\displaystyle \frac{dt^*}{dt}= \frac 12\sigma_Sm_I(t)$,  it follows that $u(x,t^*)$ satisfies 
 \be\label{heat-u}
\frac{\partial u(x,t^*)}{\partial t^*} =  \frac{\partial^2}{\partial x^2}(x^2 u(x,t^*)). 
 \ee
As shown among others in \cite{tos2016}, the solution to equation \fer{heat-u} can be easily related to the solution to the heat equation 
 \be\label{heat}
\frac{\partial v(y,t)}{\partial t} = \frac{\partial^2 v(y,t)}{\partial y^2}\,,
 \ee 
by resorting to a suitable transformation. Owing to this relationship, one shows  that the linear diffusion equation \fer{heat-u} possesses a (unique) source-type solution given by the lognormal density 
 \be\label{ln}
 L_{t^*}(x) =  \frac 1{\sqrt{4\pi t^*}\, x}\exp\left\{ - \frac{(\log\, x + t^*)^2}{4t^*} \right\},
 \ee
departing  at time $t^*=0$ from a Dirac delta function located in $x=1$.

In analogy with the heat equation \fer{heat}, where the unique solution $v(y,t)$ to the initial value problem is found to be the convolution product of the initial datum $v_0(y)$ with the source-type solution 
\be\label{max}
M_t(y) = \frac 1{\sqrt{4\pi t}}\exp\left\{ - \frac{y^2}{4t} \right\},
\ee
that is
 \[
 v(y,t) = \int_\R M_t(y-z) v_0(z) \, dz,
 \]
it is a simple exercise to verify that the unique solution to the diffusion equation \fer{heat-u} corresponding to the initial datum $u_0(x)$ is given by the expression
 \be\label{sol}
  u(x,t^*) = \int_{\R_+} \frac 1z\, u_0\left(\frac xz\right) L_{t^*}(z)  \, dz. 
 \ee
It is immediate to show that both the mass and the mean value of the solution \fer{sol} are preserved in time.
Finally, going back to the function $g_S(x,t)$ we have
\[
g_S(x,t) = \int_{\R_+} \frac 1z\, g_S\left(\frac xz, 0\right) L_{t^*}(z)  \, dz,
\]
so that
\be\label{solS}
m_S(t) f_S(m_S(t),t) = \int_{\R_+} \frac 1z\, f_S\left(\frac xz, 0\right) m_S(0) L_{t^*}(m_S(0) z)  \, dz.
\ee
Letting now $t \to \infty$ we recover the final profile of the density of susceptible in the case without re-infection. 
It is essential to remark that, in reason of the definition of $t^*$
\be\label{limS}
\lim_{t\to \infty} t^*(t) = t^*_\infty =\frac 12 \sigma_S \int_0^\infty m_I(s) ds = \frac{\sigma_S}{2\gamma} (m_R^\infty -m_R^0). 
\ee
The value of the integral in \fer{limS} follows from integration of the third equation in \fer{macro_1_red}. Hence, 
passing to the limit in \fer{solS} we obtain
\[
m_S^\infty f_S^\infty(m_S^\infty x) = \int_{\R_+} \frac 1z\, f_S\left(\frac xz, 0\right) m_S^0 L_{t^*_\infty}(m_S^0 z)  \, dz
\]
or, what is the same 
\be\label{fineS}
f_S^\infty( x) = \int_{\R_+} \frac 1z\, f_S\left(\frac xz, 0\right) \frac{m_S^0}{m_S^\infty}L_{t^*_\infty}\left(\frac{m_S^0}{m_S^\infty} z\right)  \, dz.
\ee
Hence, since $t^*_\infty$ is bounded,  the value of the asymptotic density of susceptible individuals maintains memory of the initial configuration. The passage of the epidemic wave modifies the initial density of the population according to the integral formula \fer{fineS}, which depends, in addition to the initial density of susceptible individuals,  on the initial and final values of the means of the susceptible and recovered individuals.

A similar but simpler result can be obtained for the asymptotic distribution of recovered individuals. If $\alpha = 0$, the third equation in \fer{eq:system_PDE} simplifies to 
\be\label{R-nore}
\dfrac{\partial f_R(x,t)}{\partial t} + \gamma m_I(t)\dfrac{\partial}{\partial x} \left(f_R(x,t) \right)=0\,.
\ee
Equation \fer{R-nore} can be easily solved by setting
\[
f_R(x,t) = w(x, \tau),
\]
where $ \tau = \tau (t)$ satisfies
\[
\frac{d\tau}{dt} = \gamma m_I(t).
\]
Then, $w(x,\tau)$ solves the linear transport equation
\be\label{trans}
\dfrac{\partial w(x,\tau)}{\partial \tau} + \dfrac{\partial w(x,\tau)}{\partial x}=0,
\ee
whose solution is given by $w(x,\tau)= w(x-\tau,0)$.
Consequently
\be\label{R-sol}
f_R(x,t) = f_R\left(x- \gamma\int_0^t m_I(s)\, ds , 0\right) I\left(x \ge \gamma\int_0^t m_I(s)\, ds\right),
\ee
where $I(A)$ is the indicator function of the set $A\subseteq \R_+$.
Passing to the limit in \fer{R-nore} we obtain
\be\label{R-fin}
f_R^\infty(x,t) = f_R\left(x - (m_R^\infty - m_R^0), t=0\right)I\left(x \ge m_R^\infty - m_R^0 \right).
\ee
Therefore, also the density of recovered individuals maintains memory of the initial condition.

 \begin{remark}
 One can notice that, while the density of susceptible individuals, as given by \fer{solS}, is modified in time according to the interaction with a lognormal density, the density of recovered individuals does not modify its shape, except for a time transport which clearly quantifies the increasing number of individuals which have been infected in time by the epidemic spreading.
 \end{remark} 
 
  \begin{remark}
 What obtained is consistent with classical results on the SIR model, which are translated in our setting at the level of mean quantities $m_J$, $J \in \mathcal{C}$, and for which the constant of motion is
\begin{equation}\label{constant_motion}
	\Gamma(m_S,m_I)=m_S+m_I-\dfrac{\gamma}{\beta}\log(m_S). 
\end{equation}
This quantity describes heteroclinic orbits connecting equilibria, when the basic reproduction number $\mathcal{R}_0>1$. This allows to derive a relation between $(m_S(0),m_I(0))$ and $(m_S^\infty,m_I^\infty) =
(m_S^\infty,0)$, and thus $m_R^\infty = 1-m_S^\infty$, being $m_S + m_I + m_R$  constant under the flow defined by the SIR model \eqref{macro_1_red}.  This feature is confirmed at the density level through the relation \eqref{fineS} defining the asymptotic density $f_S^\infty$, containing information not only of the mean values but on all the moments of $f_S^\infty$.
 \end{remark}

\section{Numerical tests}\label{sec:numerics}
In this section we perform several numerical tests to highlight the consistency of the proposed approach with respect to the obtained equilibria of the kinetic SIR model with reinfection. In more detail, for the evolution of the system of Fokker-Planck equations \eqref{eq:system_PDE} with no-flux boundary conditions we adopt the semi-implicit structure preserving scheme defined in \cite{PareschiZanella2018}. These methods are capable to approximate with arbitrary accuracy the large time equilibrium solution, by preserving the main physical properties in the transient regimes, like the positivity of the numerical solution. In particular, we considered a discretization of the domain $[0,100]$ obtained with $N_x = 3001$ gridpoints and $\Delta t = 0.1$. We consider a 2nd order semi-implicit scheme with a 4th order numerical approximation of the quasi-equilibrium. 

In the following tests we will fix the parameters as reported in Table \ref{tab:1}.

\begin{table}
\begin{center}
\begin{tabular}{ |c| c| c| }
\hline
Parameter & Value & Meaning \\
\hline\hline  
 $\alpha$ & 1/50 & Loss of immunity rate \\ 
 $\beta$ & 1/20 & Infection rate \\  
 $\gamma$ & 1/14 & Recovery rate   \\
  $\sigma_S$ & 1/100 & Diffusion coefficient susceptibles \\
 $\sigma_I$ & 1/100 & Diffusion coefficient infected \\
  $\sigma_R$ & 1/50 & Diffusion coefficient recovered \\
 \hline 
\end{tabular}
 \caption{Parameters chosen in the numerical tests. }
 \label{tab:1}
 \end{center}
\end{table}

We consider as initial distribution $f_J(x,0)$ the uniform distribution centered in $m_J(0)$ and with variance $\frac{1}{120}$. In particular, we consider 
\[
m_S(0) = 4,\qquad m_I(0)= 1, \qquad m_R(0) = \dfrac{1}{2}.
\]

\begin{figure}
\includegraphics[scale = 0.27]{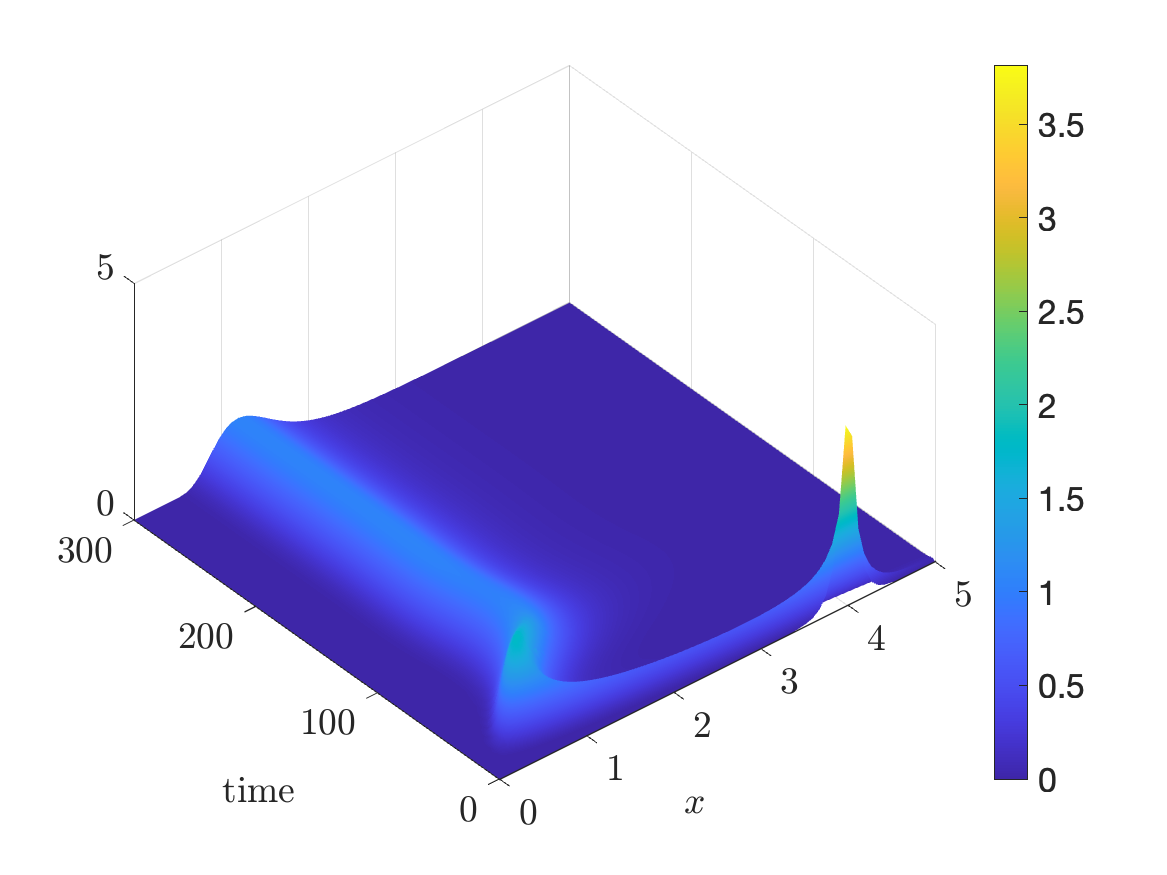}
\includegraphics[scale = 0.27]{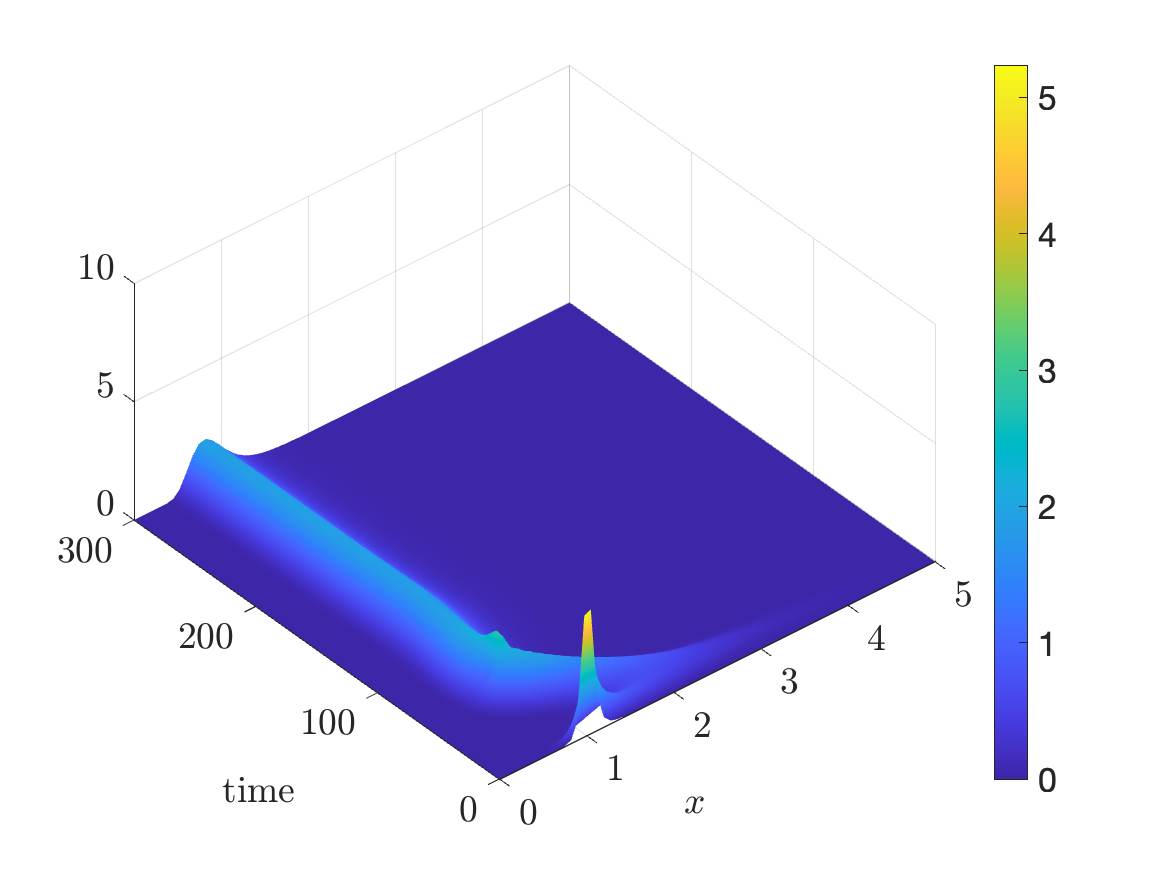}
\includegraphics[scale = 0.27]{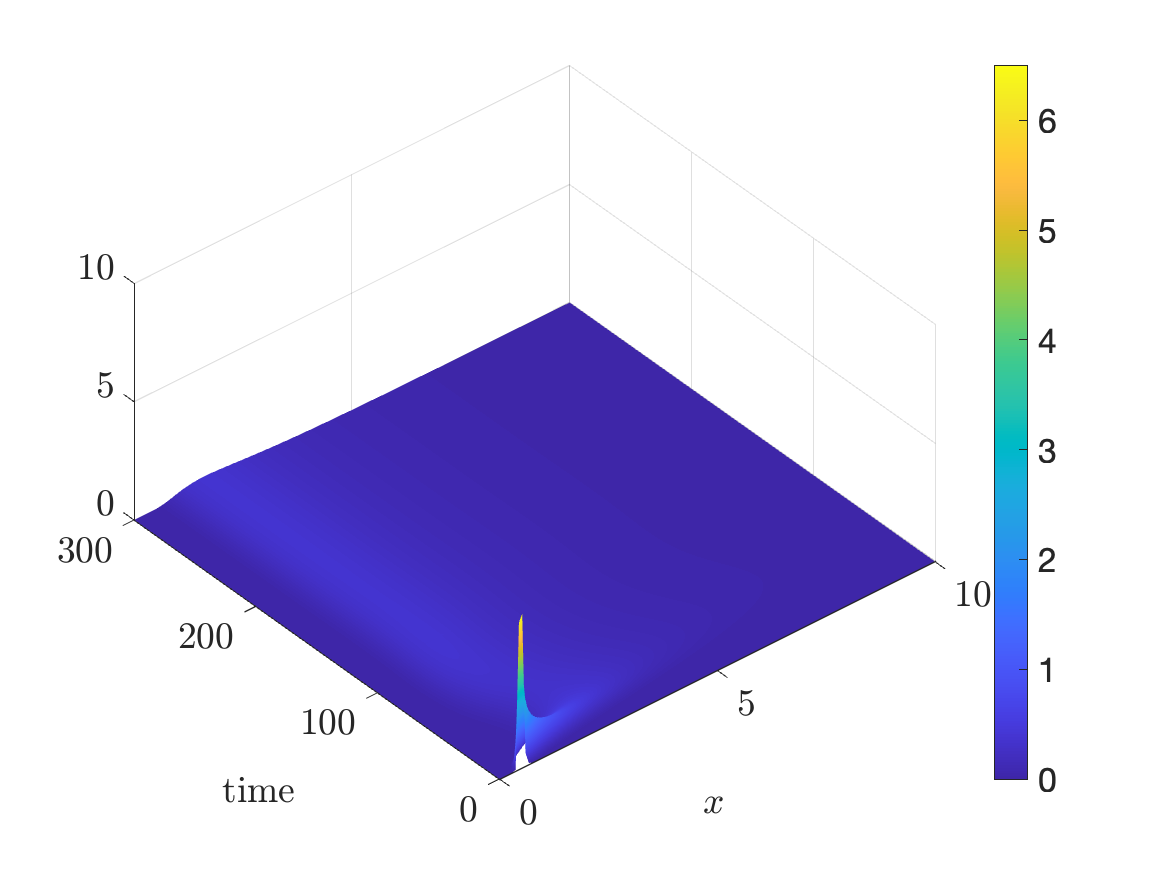}
\caption{Evolution of the kinetic distributions $f_S(x,t)$ (left), $f_I(x,t)$ (center) and $f_R(x,t)$ (right), over the time horizon $[0,100]$. }
\label{fig:1}
\end{figure}

The evolution of the kinetic distributions $f_J(x,t)$ obtained from the system \eqref{eq:system_PDE} is shown in Figure \ref{fig:1}. 

\begin{figure}
\includegraphics[scale = 0.27]{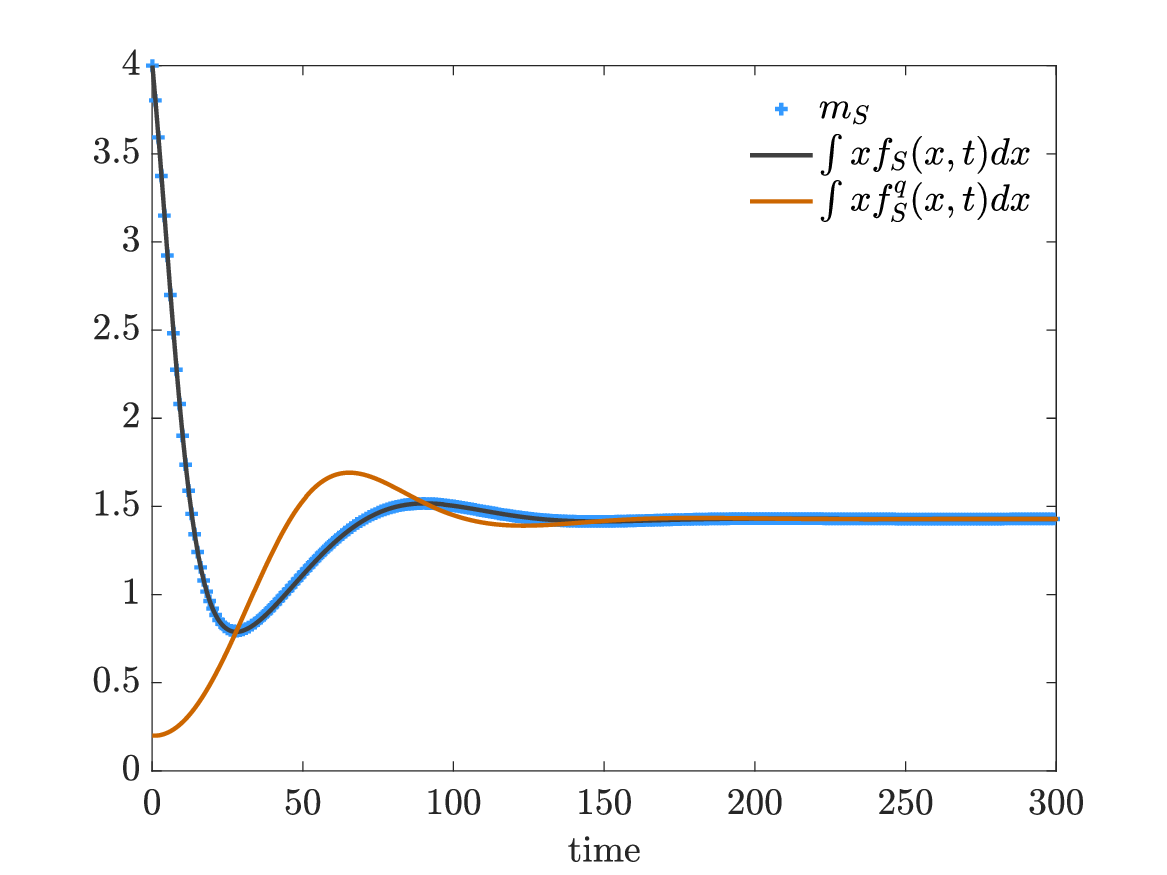}
\includegraphics[scale = 0.27]{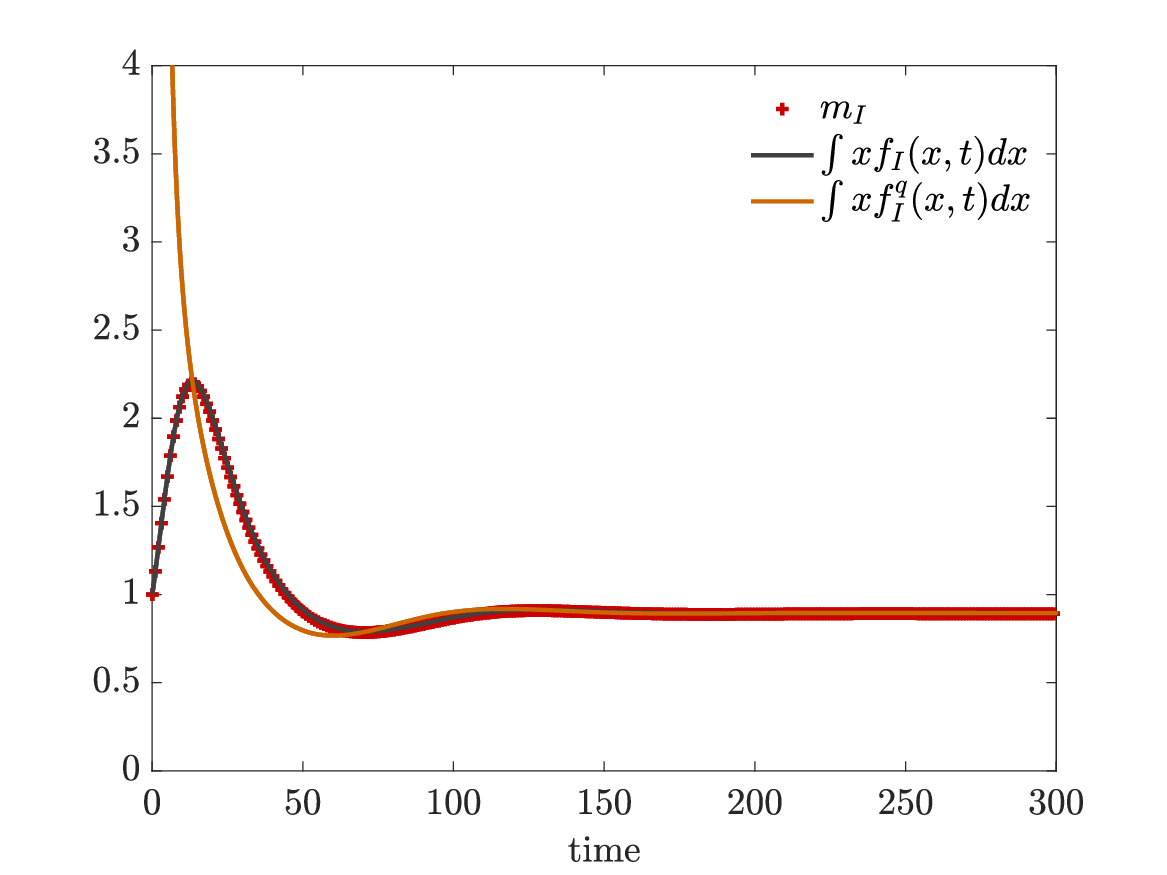}
\includegraphics[scale = 0.27]{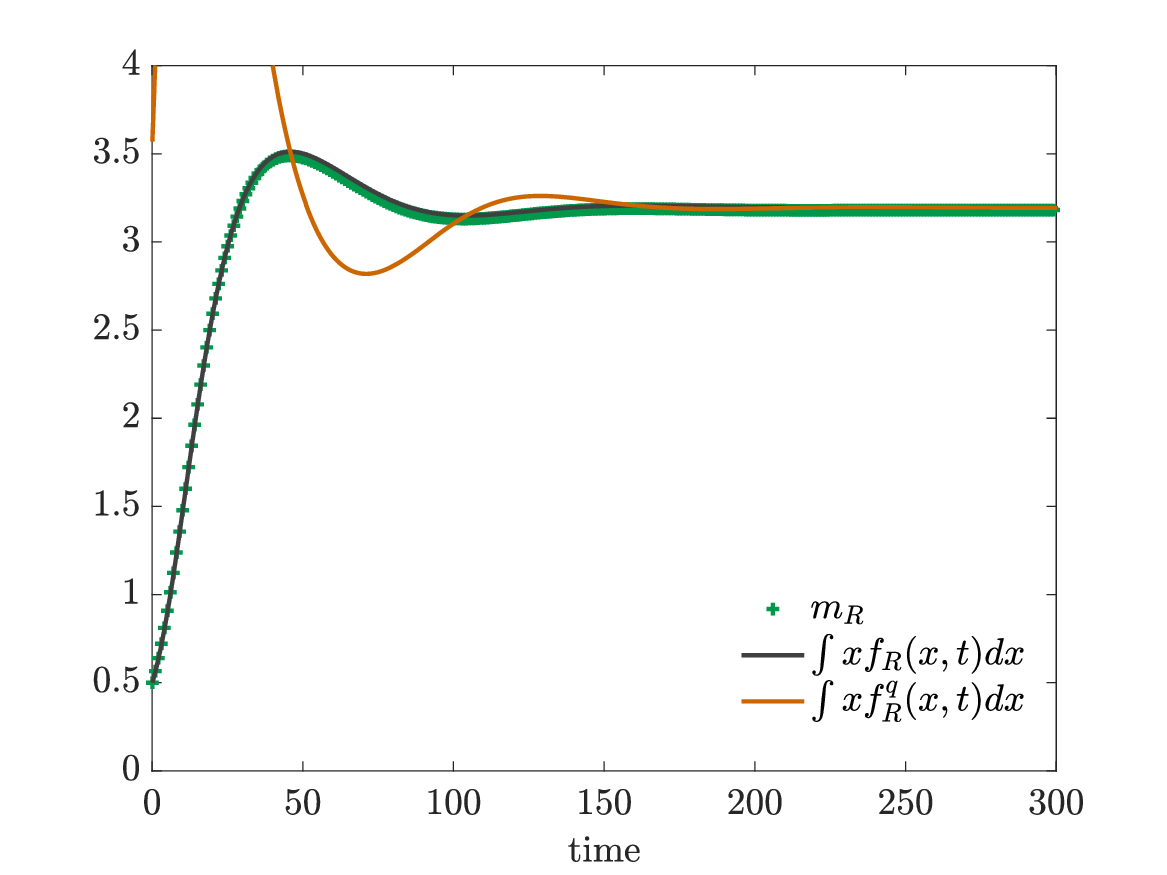}\\
\includegraphics[scale = 0.27]{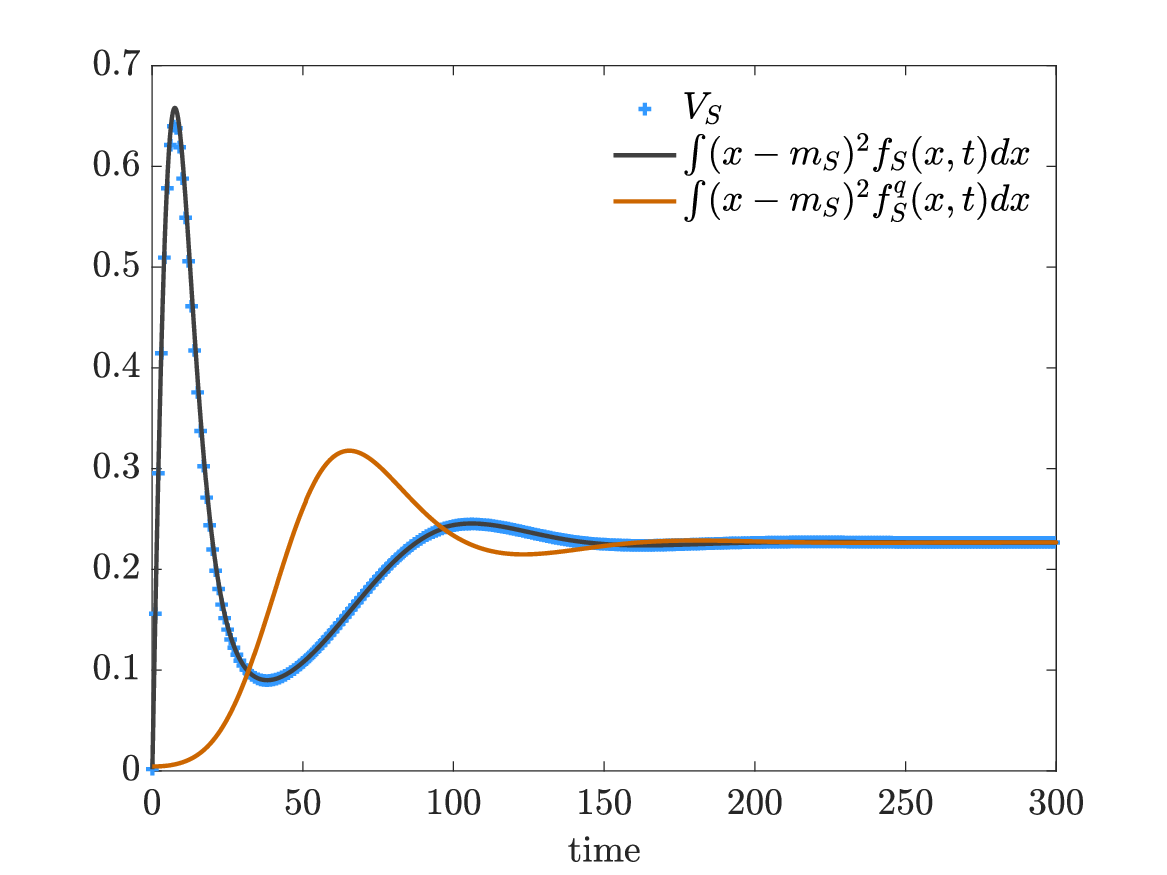}
\includegraphics[scale = 0.27]{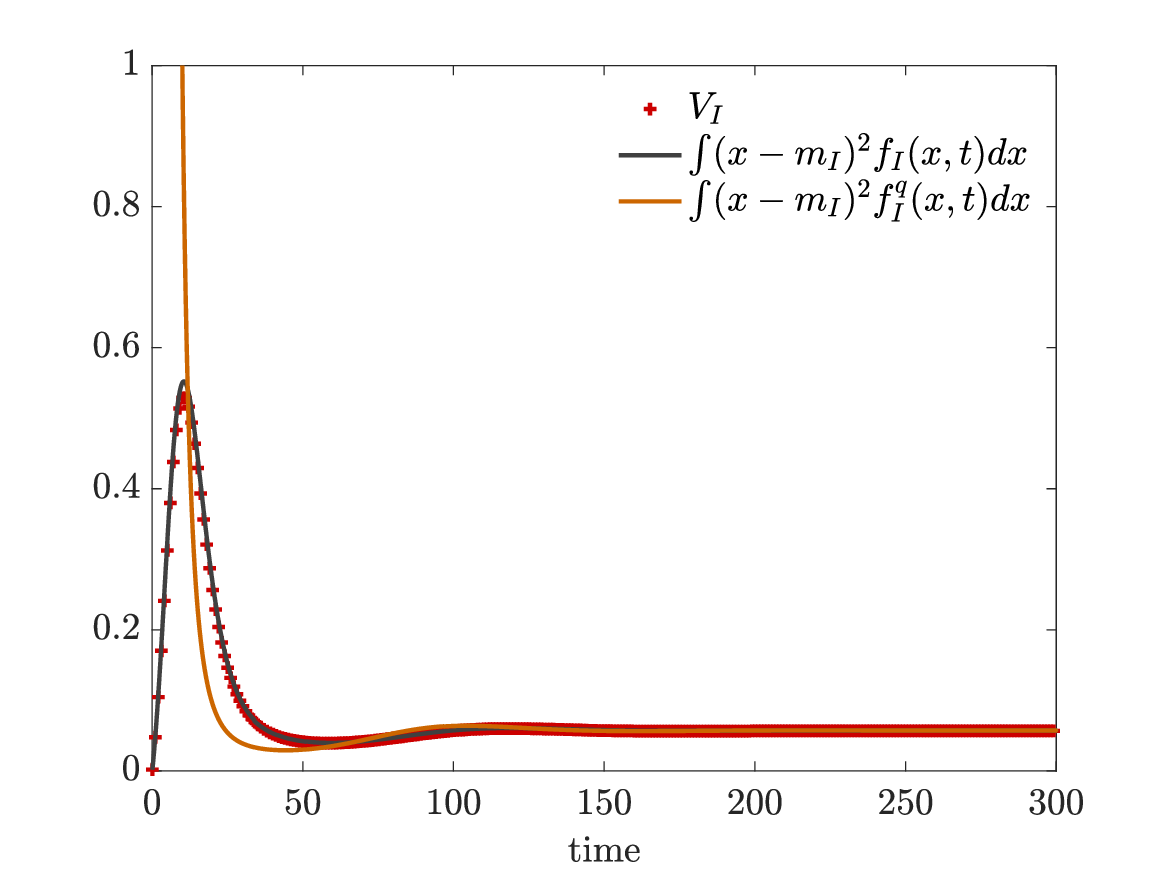}
\includegraphics[scale = 0.27]{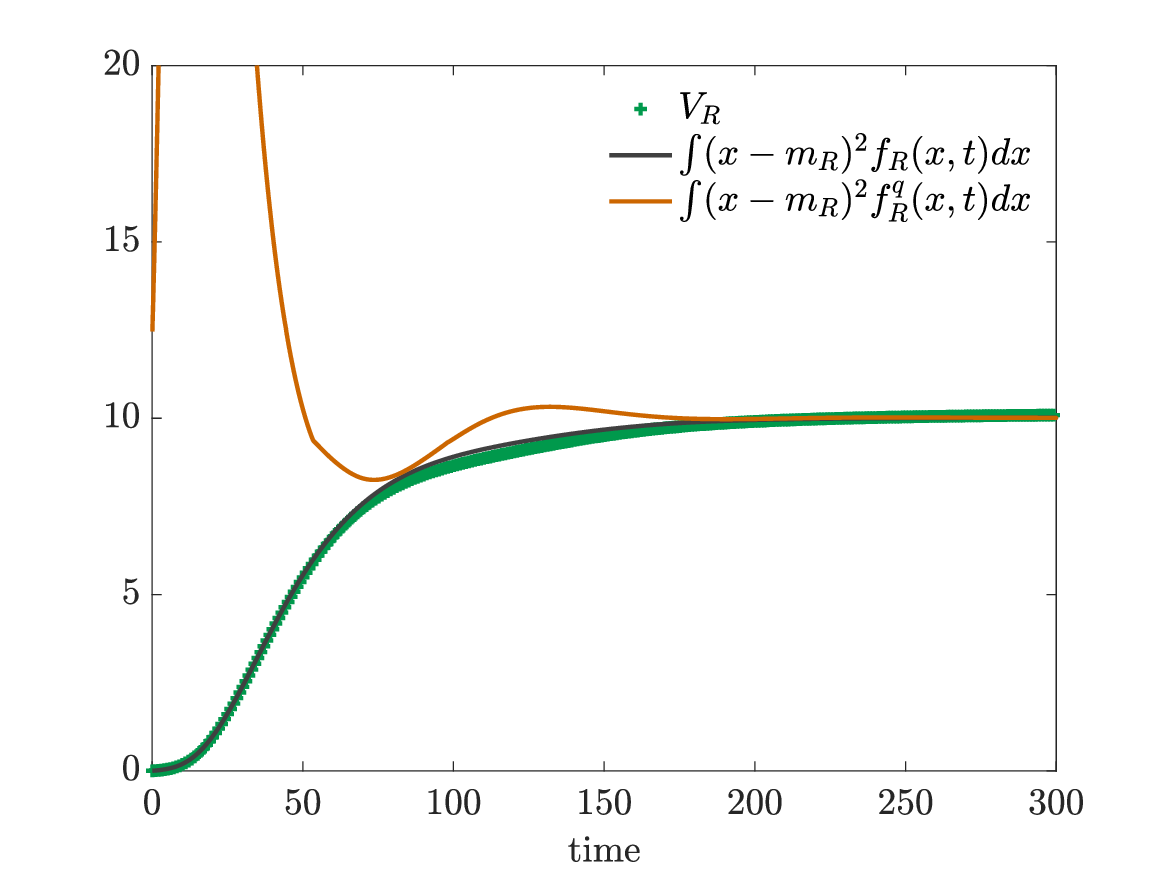}
\caption{Evolution of $m_J(t)$ and $V_J(t)$ from the solution of the system of Fokker-Planck equations and from the closed system for means \eqref{eq:mSIR} and variances \eqref{eq:variance_SIR_PDE}. We report also the evolution of the mean and variance of the quasi-equilibrium densities $f_J^q(x,t)$, $J \in \mathcal C$ to highlight the consistency of the large time behaviour.} 
\label{fig:2}
\end{figure}

In Figure \ref{fig:2} we show the evolution of mean values $m_J(t)$ solution of the SIR model \eqref{eq:mSIR} together with the mean of the approximated solution of the system of Fokker-Planck equations \eqref{eq:system_PDE} and of the quasi-equilibrium distributions \eqref{local}-\eqref{eq:localR}. We may observe good agreement between the mean and the extrapolated means of the Fokker-Planck system; the mean of the quasi-equilibrium distribution converges for large times toward the value of $m_J^\infty$ as expected. Furthermore, we compare the evolution of the obtained system for variances \eqref{eq:variance_SIR_PDE} and their approximation together with the variances of the quasi-equilibrium distributions. Their large time behaviour is consistent with the obtained equilibrium value \eqref{eq:Vinfty}. 

\begin{figure}
\includegraphics[scale = 0.27]{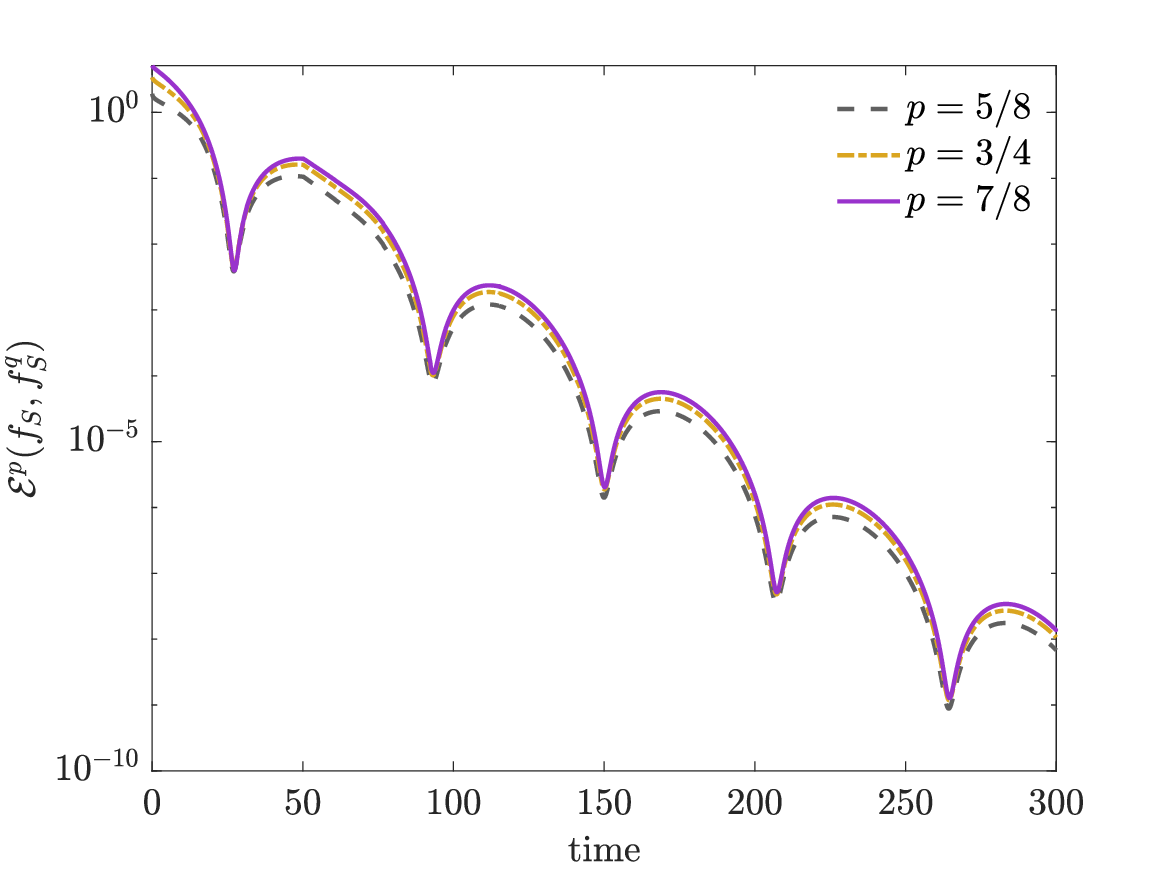}
\includegraphics[scale = 0.27]{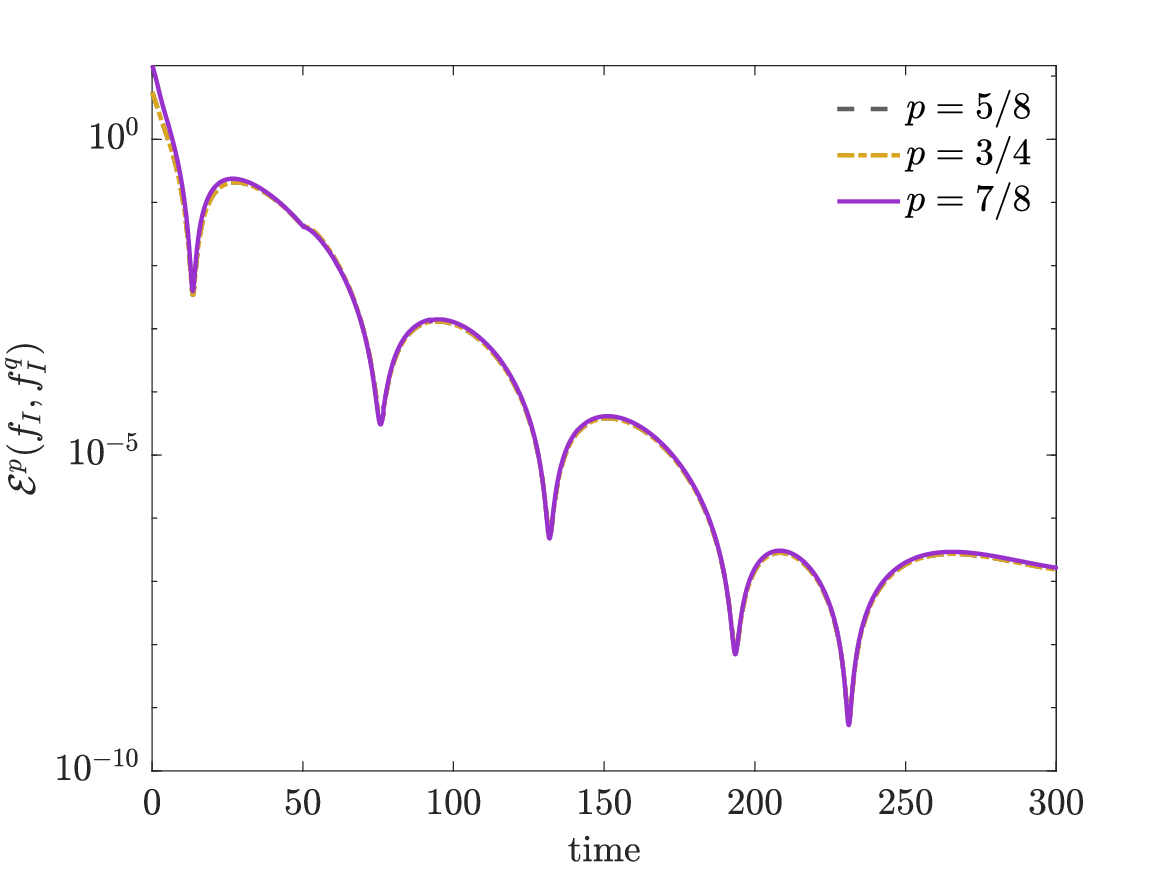}
\includegraphics[scale = 0.27]{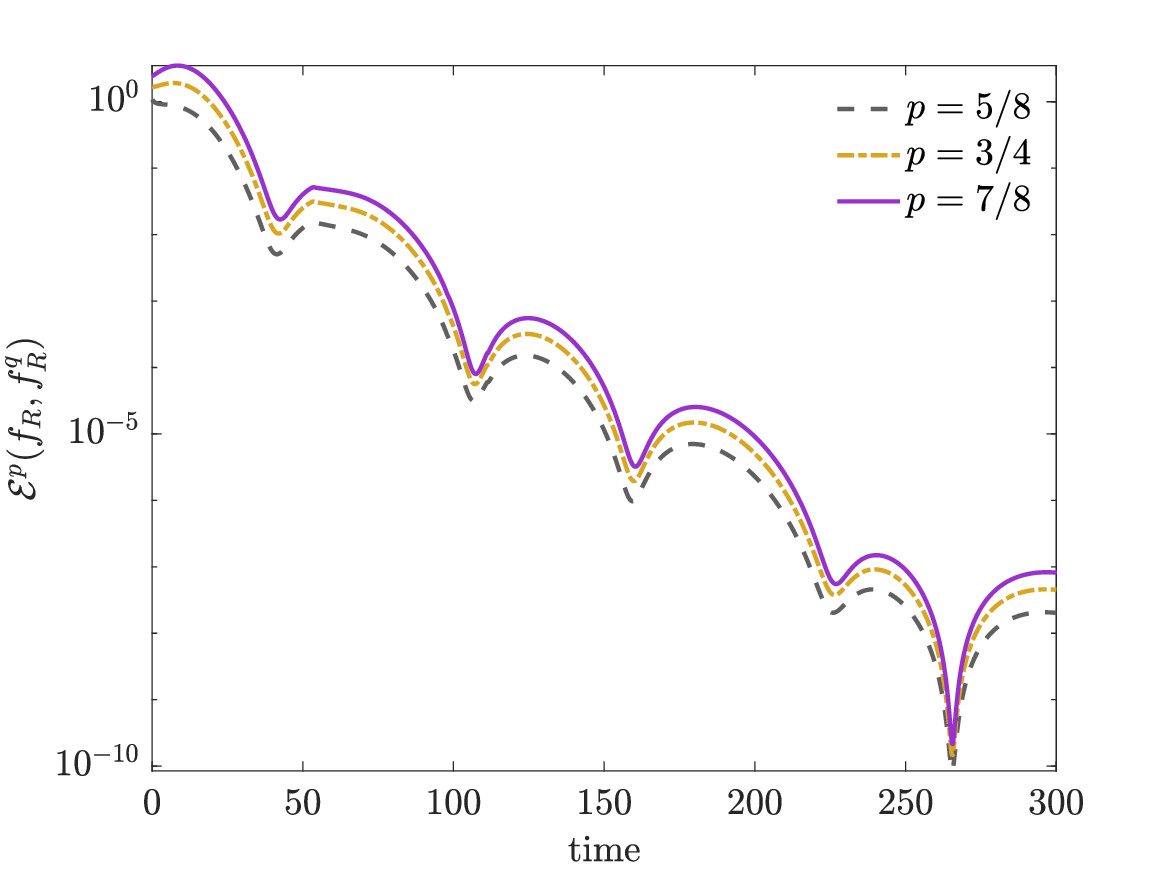}
\caption{Evolution of the energy distance $\mathcal E^p(f_J,f_J^q)$ for several weights $p = 5/8,3/4,7/8$. }
\label{fig:4}
\end{figure}

Finally, the convergence of $f_J(x,t)$ towards the quasi-equilibrium distributions $f^q_J(x,t)$ defined  in \eqref{local}-\eqref{eq:localR} is measured in terms of the energy distance $\mathcal E^p$ given in \eqref{Energy} with $p = 5/8,3/4,7/8$. In Figure \ref{fig:4} we may observe how the decay of this distance is non-monotone but it is able to provide a trend to equilibrium for the considered system of equations. 

\section*{Conclusion}

In this paper, our main goal is to understand the dynamics of multi-agent systems in the context of infectious diseases, with the aim of recovering the large-time behavior of the solutions. To this end, we introduce a kinetic framework to model the temporal evolution of the statistical distributions of the population densities in the three compartments: susceptible, infectious, and recovered, in an epidemic spread governed by susceptible-infectious interactions. The proposed model is based on a system of Boltzmann-type equations that describe binary interactions between susceptible and infectious individuals, supplemented by linear redistribution operators accounting for recovery and reinfection dynamics.
The relationship of these kinetic equations to classical compartmental epidemiological models is based on the fact that the mean values of the statistical distributions of the three populations, satisfying the kinetic system, solve a classical SIR model with or without re-infection,  where the relevant parameters are explicitly derived from the microscopic interaction rules.  

Following a classical strategy, a relevant simplification of the Boltzmann system is obtained in the grazing collision regime, where the Boltzmann system is well-approximated by a coupled system of Fokker-Planck equations. This limit allows for the analysis of the aggregate dynamics, including the large-time behavior of the population densities. In this context, we rigorously prove the convergence to equilibrium of the resulting mean-field system in a suitable Sobolev space by means of the so-called energy distance. The analysis reveals the dissipative structure of the interaction dynamics and the role of the microscopic terms in driving the system toward a stable equilibrium configuration. These results provide a multi-scale perspective connecting kinetic theory with classical epidemic models and opens new research directions both at the level of rigorous derivation of mean-field compartmental models and at the level of modelling and control of interacting many-agent systems. 

\appendix\label{app:FP}
\section{Convergence to equilibrium for Fokker--Planck equations on $\mathbb R_+$}\label{sec:wealth}

Let $f(x,t)$ denote the probability density at time $t \ge 0$ of agents with personal wealth $x\ge0$, departing from an initial density $f_0(x)$ with finite variance.
 
The evolution in time of the density $f(x,t)$ was described in \cite{bouchaud2000} by the  Fokker--Planck equation 
 \be\label{FP2c}
 \frac{\partial f(x,t)}{\partial t} =  \frac \sigma{2}\frac{\partial^2 }
{\partial x^2}\left( x^2 f(x,t)\right) +  \frac{\partial }{\partial x}\left[
(\lambda x-\mu) f(x,t)\right],
\ee
where $\sigma$, $\lambda$ and $\mu$  denote  positive constants.
The key features of equation \fer{FP2c} is  that, in presence of no flux boundary conditions at the point $x=0$, i.e.
\begin{equations}
&\left. x^2f(x,t) \right|_{x=0} = 0\,, \\
& \left. (\lambda x-\mu) f(x,t)    +  \frac \sigma{2}\frac{\partial }{\partial x}\left( x^2 f(x,t) \right)  \right|_{x=0} = 0,
 \end{equations}
the solution is mass  preserving, and approaches in time  a unique stationary solution of unit mass  \cite{torregrossa2018}. This stationary state is given by the inverse
 Gamma distribution 
 \be\label{equi2}\displaystyle
f_\infty(x) =\frac{\left( \frac {2\mu}\sigma\right)^{1+ \frac{2\lambda}\sigma }}{\Gamma\left( 1+ \frac{2\lambda}\sigma\right) }\frac 1{x^{2+ \frac{2\lambda}\sigma }} {\exp\left(-\frac{2\mu}\sigma \frac 1x \right)},
 \ee
 of shape parameter $\nu = 1+ (2\lambda)/\sigma$ and scale parameter $\omega = (2\mu)/\sigma$. This stationary state is of finite variance provided $2\lambda >\sigma$. We will always refer to a set of constants satisfying this condition, which corresponds to a suitable balance between the diffusion and the drift operators. 
 
 It is immediate to recognize that the mean value $m(t)$ evolves in time according to 
 \be\label{mean-ev}
 \frac{dm(t)}{dt} = -\lambda m(t) + \mu,
 \ee
 while the variance $V(t)$ satisifies
 \be\label{var-ev}
  \frac{dV(t)}{dt} = -(2\lambda-\sigma)  V(t) + \sigma m^2(t).
 \ee
 Clearly, $m(t) \to \mu/\lambda$ as $t \to \infty$, while, since $2\lambda >\sigma$, $V(t) \to \sigma\mu^2/ (\lambda^2(2\lambda -\sigma))$, namely the mean value and the variance of the steady state \fer{equi2}. 
 
It is interesting to observe that equation \fer{FP2c} can be expressed in terms of the Fourier transform $\fc(\xi,t)$  of $f(x,t)$, where, as usual
\be\label{fou}
\fc(\xi,t) = \int_{\R_+} f(x,t) e^{-i\xi x} \, dx.
\ee
 The Fourier transformed version of the Fokker--Planck equation \fer{FP2c} reads \cite{torregrossa2018}
\begin{equation}  \label{fourif}
\frac{\partial \widehat{f}(\xi,t)}{\partial t}= \frac{\sigma}{2}\xi^2 \dxixi \widehat{f}(\xi,t)-\lambda\xi \dxi\widehat{f}(\xi,t)-i\mu\xi \widehat{f}(\xi,t) . 
\end{equation}
Let  $ \fc(\xi,t)=a(\xi,t)+ib(\xi,t)$. Then the real and imaginary parts of $\fc$  satisfy
\begin{equations}\label{rim}
\frac{\partial a(\xi,t)}{\partial t}&=\frac{\sigma}{2}\xi^2 \dxixi a(\xi,t)-\lambda\xi \dxi a(\xi,t)+\mu\xi  b(\xi,t),  \\
\frac{\partial b(\xi,t)}{\partial t}&=\frac{\sigma}{2}\xi^2 \dxixi b(\xi,t)-\lambda\xi \dxi b(\xi,t) -\mu\xi  a(\xi,t).
\end{equations}
Let us multiply equations \fer{rim}  by $2a$ and, respectively, by $2b$. Summing up we get the evolution equation satisfied by $|\fc(\xi,t)|^2$.
\begin{equation}
\dt \fq=\sigma\xi^2\left[ a \dxixi a+b \dxixi b\right]-\lambda\xi \frac{\partial \fq}{\partial \xi}. \label{eq-f2}\\
\end{equation}
Hence, multiplying by $|\xi|^{-2p}$, $p >0$, and integrating over $\mathbb{R}$ with respect to $\xi$, we obtain the evolution equation of the $\dot H_{-p}-$norm of $f(v,t)$, where, as usual,  the homogeneous Sobolev space
 $\dot H_{-p}$, is defined by the norm
 \[
\|f\|_{\dot H_{-p}}=
\int_\R|\xi|^{-2p}\, |\hat f (\xi)|^2\,  d \xi.
 \]
We obtain
\begin{equation}
\dt \int_{\mathbb{R}}|\xi|^{-2p}  \,\fq\, d \xi=\sigma\int_{\mathbb{R}}|\xi|^{2-2p}\bigg[ a \dxixi a+b \dxixi b\bigg]d \xi-\lambda\int_{\mathbb{R}} \xi|\xi|^{-2p}  \frac{\partial \fq}{\partial \xi}d \xi, \label{eq-HP}\\
\end{equation}
and integrating by parts the two integrals, it results
\begin{equation}
\dt \int_{\mathbb{R}}|\xi|^{-2p} \,\fq\, d \xi=(1-2p)\bigg[\sigma(1-p)+\lambda\bigg]\int_{\mathbb{R}}|\xi|^{-2p}\fq d \xi-\sigma\int_{\mathbb{R}}|\xi|^{2-2p} \bigg[ \big|\dxi a\big|^2+\big|\dxi b\big|^2\bigg] d \xi
\label{eq-HP2}.
\end{equation}
Let us set $p >1/2$. Then, since for any constant $R>0$
\[
0 \le \int_{\mathbb{R}}|\xi|^{-2p}\left(\xi \frac{\partial}{\partial \xi}a -R a \right)^2  d \xi= \int_{\mathbb{R}}|\xi|^{2-2p}  \big|\dxi a\big|^2 +\left[R^2+ R(1-2p)\right] \int_\R |\xi|^{-2p} a^2(\xi) \, d\xi,
\]
it follows that
\[
\int_{\mathbb{R}}|\xi|^{2-2p}  \big|\dxi a\big|^2  \ge \left[ R(2p-1)- R^2\right] \int_\R |\xi|^{-2p} a^2(\xi) \, d\xi.
\]
Same inequality for $b(\xi)$. Optimizing over $R$ we finally obtain  \cite{torregrossa2018}
$$
\int_{\mathbb{R}}|\xi|^{2-2p} \bigg[ \big|\dxi a\big|^2+\big|\dxi b\big|^2\bigg] d \xi\geq \frac{(2p-1)^2}{4}\int_{\mathbb{R}}|\xi|^{-2p}\fq\, d\xi.
$$
Therefore
 \be
\dt \int_{\mathbb{R}} |\xi|^{-2p} \fq\, d \xi \leq -(2p-1)\bigg[ \sigma\frac{3-2p}4 + \lambda \bigg]\int_{\mathbb{R}}|\xi|^{-2p} \fq d \xi.
\label{eq-HP3}
\end{equation}
Hence, provided $1/2 <p\le 3/2$
inequality \eqref{eq-HP3} implies that if the initial data has bounded $\dot H_{-p}-$norm, then for all $t>0$, the $\dot H_{-p}-$norm of the solution remains uniformly bounded. The previous computations can be fruitfully used to show that the solution to equation \fer{FP2c} converges exponentially towards the steady state in the Sobolev space $H_{-p}$.   \cite{torregrossa2018}. 
Since the Fokker--Planck equation \fer{FP2c} is linear, we can repeat the previous computations for the difference $\hat f(\xi,t) - \hat f_\infty(\xi)$, which clearly solve equation \fer{fourif}. Consequently, thanks to \fer{eq-HP3}, if initially the  $H_{-p}$ norm of the difference $\hat f(\xi,t) - \hat f_\infty(\xi)$ is initially bounded for some $p$ in the interval $1/2 < p \le 3/2$,  this norm is exponentially decaying to zero as $t \to \infty$.

As shown in \cite{szekely2013}, if $1/2 < p< 3/2$, the  $H_{-p}$ norm can be equivalently expressed in the physical space according to the formula
\begin{equations}\label{Energy}
\mathcal{E}^p(f,f_\infty) &= 2 \int_{\R_+\times\R_+} |x-y|^{2p-1} f(x)f_\infty(y) \, dxdy - \\
&  \int_{\R_+\times\R_+} |x-y|^{2p-1} f(x)f(y) \, dxdy - \int_{\R_+\times\R_+} |x-y|^{2p-1} f_\infty(x)f_\infty(y) \, dxdy = \\
& \frac{2p-1}{\sqrt 2} \frac{2^{2p-1} \Gamma\left( p\right)}{ \Gamma\left( \frac 32- p\right)}\int_\R \frac{|\hat f(\xi,t) - \hat f_\infty(\xi)|^2}{|\xi|^{2p}}.
\end{equations}
In statistics, the expression $\mathcal{E}^p(f,f_\infty)$ is known under the name of \emph{energy distance}, and allows for numerical computations without passing to Fourier transform. 
\medskip

\bigskip
\emph {Acknowledgments} 
This work has been written within the activities of GNFM group of INdAM (National Institute of High Mathematics). M.Z. acknowledges partial support by PRIN2022PNRR project No.P2022Z7ZAJ, European Union - NextGenerationEU and by ICSC - Centro Nazionale di Ricerca in High Performance Computing, Big Data and Quantum Computing, funded by European Union - NextGenerationEU.

\setlength\parindent{0pt}
\bibliographystyle{plain}
\bibliography{biblio.bib}

@book {MR1307620,
    AUTHOR = {Cercignani, C. and Illner, R. and Pulvirenti, M.},
     TITLE = {The mathematical theory of dilute gases},
    SERIES = {Applied Mathematical Sciences},
    VOLUME = {106},
 PUBLISHER = {Springer-Verlag, New York},
      YEAR = {1994},
     PAGES = {viii+347},
      ISBN = {0-387-94294-7},
   MRCLASS = {82C40 (76-02 76P05 82-02 82B40)},
  MRNUMBER = {1307620},
MRREVIEWER = {Giuseppe\ Toscani},
       DOI = {10.1007/978-1-4419-8524-8},
       URL = {https://doi.org/10.1007/978-1-4419-8524-8},
}

@article{bernardi2025heterogeneously,
	title={Heterogeneously Structured Compartmental Models of Epidemiological Systems: From Individual-Level Processes to Population-Scale Dynamics},
	author={Bernardi, E. and Lorenzi, T. and Sensi, M. and Tosin, A.},
	journal={Studies in Applied Mathematics},
	volume={155},
	number={2},
	pages={e70091},
	year={2025},
	publisher={Wiley Online Library}
}

@article{preziosi,
	author = {Preziosi, L. and Toscani, G. and Zanella, M.},
	doi = {https://doi.org/10.1016/j.jtbi.2021.110579},
	issn = {0022-5193},
	journal = {J. Theoret. Biol.},
	pages = {110579},
	title = {Control of tumor growth distributions through kinetic methods},
	url = {https://www.sciencedirect.com/science/article/pii/S0022519321000011},
	volume = {514},
	year = {2021}}

@article{bertaglia24,
	author = {Bertaglia, G. and Bondesan, A. and Burini, D. and Eftimie, R. and Pareschi, L. and Toscani, G.},
	doi = {10.1142/S0218202524500301},
	eprint = {https://doi.org/10.1142/S0218202524500301},
	journal = {Math. Mod. Meth Appl. Scie.},
	number = {11},
	pages = {1995-2054},
	title = {New trends on the systems approach to modeling {SARS}-{C}o{V}-2 pandemics in a globally connected planet},
	url = {https://doi.org/10.1142/S0218202524500301},
	volume = {34},
	year = {2024}}

@article{ha22,
	author = {Ha, S.-Y. and Park, H. and Yang, S.},
	doi = {10.3934/cpaa.2022127},
	fjournal = {Communications on Pure and Applied Analysis},
	issn = {1534-0392,1553-5258},
	journal = {Commun. Pure Appl. Anal.},
	mrclass = {92D30 (34D06 70G60)},
	mrnumber = {4512995},
	mrreviewer = {Yilun\ Shang},
	number = {11},
	pages = {3887--3918},
	title = {Relaxation dynamics of {SIR}-flocks with random epidemic states},
	url = {https://doi.org/10.3934/cpaa.2022127},
	volume = {21},
	year = {2022}}

@article{MR3067586,
	author = {Degond, P. and Frouvelle, A. and Liu, J.-G.},
	doi = {10.1007/s00332-012-9157-y},
	fjournal = {Journal of Nonlinear Science},
	issn = {0938-8974,1432-1467},
	journal = {J. Nonlinear Sci.},
	mrclass = {82C26 (35B25 35K55 35L60 82C22 92D50)},
	mrnumber = {3067586},
	mrreviewer = {Fernando\ Pestana\ da Costa},
	number = {3},
	pages = {427--456},
	title = {Macroscopic limits and phase transition in a system of self-propelled particles},
	url = {https://doi.org/10.1007/s00332-012-9157-y},
	volume = {23},
	year = {2013}}

@article{MR4251319,
	author = {Burger, M. and Kreusser, L. M. and Totzeck, C.},
	doi = {10.1051/cocv/2021034},
	fjournal = {ESAIM. Control, Optimisation and Calculus of Variations},
	issn = {1292-8119,1262-3377},
	journal = {ESAIM Control Optim. Calc. Var.},
	mrclass = {49N80 (70F10 82C22 92C15)},
	mrnumber = {4251319},
	mrreviewer = {Laurent\ Pfeiffer},
	pages = {Paper No. 40, 24},
	title = {Mean-field optimal control for biological pattern formation},
	url = {https://doi.org/10.1051/cocv/2021034},
	volume = {27},
	year = {2021}}

@article{during09,
	author = {D\"uring, B. and Markowich, P. and Pietschmann, J.-F. and Wolfram, M.-T.},
	doi = {10.1098/rspa.2009.0239},
	fjournal = {Proceedings of The Royal Society of London. Series A. Mathematical, Physical and Engineering Sciences},
	issn = {1364-5021,1471-2946},
	journal = {Proc. R. Soc. Lond. Ser. A Math. Phys. Eng. Sci.},
	mrclass = {91B80 (82C40 91D10)},
	mrnumber = {2552289},
	number = {2112},
	pages = {3687--3708},
	title = {Boltzmann and {F}okker-{P}lanck equations modelling opinion formation in the presence of strong leaders},
	url = {https://doi.org/10.1098/rspa.2009.0239},
	volume = {465},
	year = {2009}}

@article{during15,
	author = {D\"uring, B. and Wolfram, M.-T.},
	doi = {10.1098/rspa.2015.0345},
	fjournal = {Proceedings A},
	issn = {1364-5021,1471-2946},
	journal = {Proc. A.},
	mrclass = {91D10 (35Q20 35Q91)},
	mrnumber = {3420842},
	number = {2182},
	pages = {20150345, 21},
	title = {Opinion dynamics: inhomogeneous {B}oltzmann-type equations modelling opinion leadership and political segregation},
	url = {https://doi.org/10.1098/rspa.2015.0345},
	volume = {471},
	year = {2015}}

@article{dellamarca23,
	author = {Della Marca, R. and Loy, N. and Tosin, A.},
	doi = {10.1007/s00285-023-01901-z},
	fjournal = {Journal of Mathematical Biology},
	issn = {0303-6812,1432-1416},
	journal = {J. Math. Biol.},
	mrclass = {92D30 (35Q20 35Q70 35Q84 37N25)},
	mrnumber = {4568210},
	mrreviewer = {Rafael\ Bravo de la Parra},
	number = {4},
	pages = {Paper No. 61, 28},
	title = {An {SIR} model with viral load-dependent transmission},
	url = {https://doi.org/10.1007/s00285-023-01901-z},
	volume = {86},
	year = {2023}}

@article{lorenzi24,
	author = {Lorenzi, T. and Paparelli, E. and Tosin, A.},
	doi = {10.4310/cms.240918045626},
	fjournal = {Communications in Mathematical Sciences},
	issn = {1539-6746,1945-0796},
	journal = {Commun. Math. Sci.},
	mrclass = {92D30 (35Q92 35R09 45K05)},
	mrnumber = {4836954},
	mrreviewer = {Ran\ Zhang},
	number = {8},
	pages = {2131--2165},
	title = {Modelling coevolutionary dynamics in heterogeneous {SI} epidemiological systems across scales},
	url = {https://doi.org/10.4310/cms.240918045626},
	volume = {22},
	year = {2024}}

@article{carrillo10,
	author = {Carrillo, J. A. and Fornasier, M. and Rosado, J. and Toscani, G.},
	doi = {10.1137/090757290},
	fjournal = {SIAM Journal on Mathematical Analysis},
	issn = {0036-1410,1095-7154},
	journal = {SIAM J. Math. Anal.},
	mrclass = {35B40 (35Q92 82D99 92D50)},
	mrnumber = {2596552},
	mrreviewer = {Giuseppe\ Maria\ Coclite},
	number = {1},
	pages = {218--236},
	title = {Asymptotic flocking dynamics for the kinetic {C}ucker-{S}male model},
	url = {https://doi.org/10.1137/090757290},
	volume = {42},
	year = {2010}}

@article{motsch14,
	author = {Motsch, S. and Tadmor, E.},
	doi = {10.1137/120901866},
	fjournal = {SIAM Review},
	issn = {0036-1445,1095-7200},
	journal = {SIAM Rev.},
	mrclass = {92D25 (74A25)},
	mrnumber = {3274797},
	number = {4},
	pages = {577--621},
	title = {Heterophilious dynamics enhances consensus},
	url = {https://doi.org/10.1137/120901866},
	volume = {56},
	year = {2014}}

@article{degond08,
	author = {Degond, P. and Motsch, S.},
	doi = {10.1142/S0218202508003005},
	fjournal = {Mathematical Models and Methods in Applied Sciences},
	issn = {0218-2025,1793-6314},
	journal = {Math. Models Methods Appl. Sci.},
	mrclass = {35Q80 (35B40 82C22 82C70 92D50)},
	mrnumber = {2438213},
	mrreviewer = {Sergey\ A.\ Vakulenko},
	pages = {1193--1215},
	title = {Continuum limit of self-driven particles with orientation interaction},
	url = {https://doi.org/10.1142/S0218202508003005},
	volume = {18},
	year = {2008}}

@article{ha09,
	author = {Ha, S.-Y. and Liu, J.-G.},
	doi = {10.4310/cms.2009.v7.n2.a2},
	fjournal = {Communications in Mathematical Sciences},
	issn = {1539-6746,1945-0796},
	journal = {Commun. Math. Sci.},
	mrclass = {82C22 (92C15 92D50)},
	mrnumber = {2536440},
	number = {2},
	pages = {297--325},
	title = {A simple proof of the {C}ucker-{S}male flocking dynamics and mean-field limit},
	url = {https://doi.org/10.4310/cms.2009.v7.n2.a2},
	volume = {7},
	year = {2009}}

@article{bolley12,
	author = {Bolley, F. and Ca\~nizo, J. A. and Carrillo, J. A.},
	doi = {10.1016/j.aml.2011.09.011},
	fjournal = {Applied Mathematics Letters. An International Journal of Rapid Publication},
	issn = {0893-9659,1873-5452},
	journal = {Appl. Math. Lett.},
	mrclass = {60K35 (35A01 35A02 35R60 60H10)},
	mrnumber = {2855983},
	number = {3},
	pages = {339--343},
	title = {Mean-field limit for the stochastic {V}icsek model},
	url = {https://doi.org/10.1016/j.aml.2011.09.011},
	volume = {25},
	year = {2012}}

@article{Toscani:1999aa,
	author = {Toscani, G. and Villani, C.},
	date = {1999/06/01},
	doi = {10.1007/s002200050631},
	id = {Toscani1999},
	isbn = {1432-0916},
	journal = {Commun. Math. Phys.},
	number = {3},
	pages = {667--706},
	title = {Sharp Entropy Dissipation Bounds and Explicit Rate of Trend to Equilibrium for the Spatially Homogeneous Boltzmann Equation},
	url = {https://doi.org/10.1007/s002200050631},
	volume = {203},
	year = {1999}}

@article{Furioli,
	author = {Furioli, G. and Pulvirenti, A. and Terraneo, E. and Toscani, G.},
	doi = {10.1142/S0218202517400048},
	eprint = {https://doi.org/10.1142/S0218202517400048},
	journal = {Math. Mod. Meth. Appl. Scie.},
	number = {01},
	pages = {115-158},
	title = {{F}okker-{P}lanck equations in the modeling of socio-economic phenomena},
	url = {https://doi.org/10.1142/S0218202517400048},
	volume = {27},
	year = {2017}}

@article{Toscani99,
	author = {Toscani, G.},
	doi = {10.1090/qam/1704435},
	journal = {Quart. Appl. Math.},
	number = {3},
	pages = {521--541},
	title = {Entropy production and the rate of convergence to equilibrium for the {F}okker-{P}lanck equation},
	volume = {LVII},
	year = {1999}}

@incollection{DeJD,
	author = {De Jong, M. C. and Diekmann, O. and Heesterbeek, H.},
	booktitle = {Epidemic Models: Their Structure and Relation to Data},
	editor = {Mollison, D.},
	pages = {84--94},
	title = {How does transmission of infection depend on population size},
	year = {1995}}

@article{riley07,
	author = {Riley, S.},
	doi = {10.1126/science.1134695},
	eprint = {https://www.science.org/doi/pdf/10.1126/science.1134695},
	journal = {Science},
	number = {5829},
	pages = {1298-1301},
	title = {Large-Scale Spatial-Transmission Models of Infectious Disease},
	url = {https://www.science.org/doi/abs/10.1126/science.1134695},
	volume = {316},
	year = {2007}}

@article{HU2013125,
	author = {Hu, H. and Nigmatulina, K. and Eckhoff, P.},
	doi = {10.1016/j.mbs.2013.04.013},
	issn = {0025-5564},
	journal = {Math. Biosci.},
	number = {2},
	pages = {125-134},
	title = {The scaling of contact rates with population density for the infectious disease models},
	url = {https://www.sciencedirect.com/science/article/pii/S0025556413001235},
	volume = {244},
	year = {2013}}

@article{Difrancesco25,
	author = {Di Francesco, M. and Zefreh, F. G.},
	doi = {10.3934/dcds.2025054},
	fjournal = {Discrete and Continuous Dynamical Systems. Series A},
	issn = {1078-0947,1553-5231},
	journal = {Discrete Contin. Dyn. Syst.},
	mrclass = {35F55 (35A01 35B40 35K40 45K05 92D30)},
	mrnumber = {4913474},
	number = {11},
	pages = {4202--4237},
	title = {Kermack-{M}c{K}endrick type models for epidemics with nonlocal aggregation terms},
	url = {https://doi.org/10.3934/dcds.2025054},
	volume = {45},
	year = {2025}}

@article{BisiLorenzani,
	author = {Bisi, M. and Lorenzani, S.},
	doi = {10.1007/s00332-024-10062-2},
	fjournal = {Journal of Nonlinear Science},
	issn = {0938-8974,1432-1467},
	journal = {J. Nonlinear Sci.},
	mrclass = {35Q20 (82C40 92D30)},
	mrnumber = {4770477},
	number = {5},
	pages = {Paper No. 84, 44},
	title = {Mathematical models for the large spread of a contact-based infection: a statistical mechanics approach},
	url = {https://doi.org/10.1007/s00332-024-10062-2},
	volume = {34},
	year = {2024}}

@book{iannellipugliese,
	author = {Iannelli, M. and Pugliese, A.},
	doi = {10.1007/978-3-319-03026-5},
	publisher = {Springer Cham},
	title = {An Introduction to Mathematical Population Dynamics: Along the trail of Volterra and Lotka},
	year = {2014}}

@article{PareschiZanella2018,
	author = {Pareschi, L. and Zanella, M.},
	date = {2018/03/01},
	doi = {10.1007/s10915-017-0510-z},
	id = {Pareschi2018},
	isbn = {1573-7691},
	journal = {J. Sci. Comput.},
	number = {3},
	pages = {1575--1600},
	title = {Structure Preserving Schemes for Nonlinear {F}okker--{P}lanck Equations and Applications},
	url = {https://doi.org/10.1007/s10915-017-0510-z},
	volume = {74},
	year = {2018}}

@article{TosZan,
	author = {Toscani, G. and Zanella, M.},
	journal = {Riv. Mat. Univ. Parma},
	number = {1},
	pages = {61--77},
	title = {On a kinetic description of {L}otka--{V}olterra dynamics},
	volume = {15},
	year = {2024}}

@article{BonMenTosZan,
	author = {Bondesan, A. and Menale, M. and Toscani, G. and Zanella, M.},
	doi = {10.1088/1361-6544/addfa1},
	journal = {Nonlinearity},
	number = {7},
	pages = {075026},
	title = {{L}otka--{V}olterra-type kinetic equations for competing species},
	volume = {38},
	year = {2025}}

@article{holling1959components,
	author = {Holling, C.S.},
	doi = {10.4039/Ent91293-5},
	journal = {Can. Entomol.},
	number = {5},
	pages = {293--320},
	publisher = {Cambridge University Press},
	title = {The components of predation as revealed by a study of small-mammal predation of the {E}uropean {P}ine {S}awfly1},
	volume = {91},
	year = {1959}}

@article{holling1966functional,
	author = {Holling, C.S.},
	doi = {10.4039/entm9848fv},
	journal = {The Memoirs of the Entomological Society of Canada},
	number = {S48},
	pages = {5--86},
	publisher = {Cambridge University Press},
	title = {The functional response of invertebrate predators to prey density},
	volume = {98},
	year = {1966}}

@article{beddington1975mutual,
	author = {Beddington, J.R.},
	doi = {10.2307/3866},
	journal = {J. Anim. Ecol.},
	pages = {331--340},
	publisher = {JSTOR},
	title = {Mutual interference between parasites or predators and its effect on searching efficiency},
	year = {1975}}

@incollection{kuno1987principles,
	author = {Kuno, E.},
	booktitle = {Advances in Ecological Research},
	doi = {10.1016/S0065-2504(08)60090-2},
	pages = {249--337},
	publisher = {Elsevier},
	title = {Principles of predator--prey interaction in theoretical, experimental, and natural population systems},
	volume = {16},
	year = {1987}}

@article{abrams2000nature,
	author = {Abrams, P.A. and Ginzburg, L.R.},
	doi = {10.1016/S0169-5347(00)01908-X},
	journal = {Trends Ecol. Evol.},
	number = {8},
	pages = {337--341},
	publisher = {Elsevier},
	title = {The nature of predation: prey dependent, ratio dependent or neither?},
	volume = {15},
	year = {2000}}

@book{ParTos-2013,
	author = {Pareschi, L. and Toscani, G.},
	publisher = {OUP Oxford},
	title = {Interacting {M}ultiagent {S}ystems: {K}inetic {E}quations and {M}onte {C}arlo {M}ethods},
	year = {2013}}

@article{dimarco2021kinetic,
	author = {Dimarco, G. and Perthame, B. and Toscani, G. and Zanella, M.},
	doi = {10.1007/s00285-021-01630-1},
	journal = {J. Math. Biol.},
	number = {1},
	pages = {4},
	publisher = {Springer},
	title = {Kinetic models for epidemic dynamics with social heterogeneity},
	volume = {83},
	year = {2021}}

@incollection{albi2022kinetic,
	author = {Albi, G. and Bertaglia, G. and Boscheri, W. and Dimarco, G. and Pareschi, L. and Toscani, G. and Zanella, M.},
	booktitle = {Predicting Pandemics in a Globally Connected World, Volume 1: Toward a Multiscale, Multidisciplinary Framework through Modeling and Simulation},
	doi = {10.1007/978-3-030-96562-4_3},
	pages = {43--108},
	publisher = {Springer},
	title = {Kinetic modelling of epidemic dynamics: social contacts, control with uncertain data, and multiscale spatial dynamics},
	year = {2022}}

@article{kermack1927,
	author = {Kermack, W.O. and McKendrick, A.G.},
	doi = {10.1098/rspa.1927.0118},
	journal = {Proc. R. Soc. London A},
	number = {772},
	pages = {700--721},
	publisher = {The Royal Society London},
	title = {A contribution to the mathematical theory of epidemics},
	volume = {115},
	year = {1927}}

@article{zanella2023,
	author = {Zanella, M.},
	doi = {10.1007/s11538-023-01147-2},
	journal = {Bull. Math. Biol.},
	number = {5},
	pages = {36},
	publisher = {Springer},
	title = {Kinetic models for epidemic dynamics in the presence of opinion polarization},
	volume = {85},
	year = {2023}}

@article{hethcote2000,
	author = {Hethcote, H.W.},
	doi = {10.1137/S0036144500371907},
	journal = {SIAM Rev.},
	number = {4},
	pages = {599--653},
	publisher = {SIAM},
	title = {The mathematics of infectious diseases},
	volume = {42},
	year = {2000}}

@article{bisi2009,
	author = {Bisi, M. and Spiga, G. and Toscani, G.},
	journal = {Commun. Math. Sci.},
	number = {4},
	pages = {901--916},
	title = {Kinetic models of conservative economies with wealth redistribution},
	volume = {7},
	year = {2009}}

@article{bouchaud2000,
	author = {Bouchaud, J-P. and M{\'e}zard, M.},
	doi = {10.1016/S0378-4371(00)00205-3},
	journal = {Phys. A},
	number = {3-4},
	pages = {536--545},
	publisher = {Elsevier},
	title = {Wealth condensation in a simple model of economy},
	volume = {282},
	year = {2000}}

@article{cordier2005,
	author = {Cordier, S. and Pareschi, L. and Toscani, G.},
	doi = {10.1007/s10955-005-5456-0},
	journal = {J. Stat. Phys.},
	pages = {253--277},
	publisher = {Springer},
	title = {On a kinetic model for a simple market economy},
	volume = {120},
	year = {2005}}

@article{torregrossa2018,
	author = {Torregrossa, M. and Toscani, G.},
	journal = {Commun. Math. Sci.},
	number = {2},
	pages = {537--560},
	title = {Wealth distribution in presence of debts. {A} {F}okker--{P}lanck description},
	volume = {16},
	year = {2018}}

@article{tos2016,
	author = {Toscani, G.},
	doi = {10.1016/j.physa.2016.06.063},
	journal = {Phys. A},
	pages = {802--811},
	publisher = {Elsevier},
	title = {Kinetic and mean field description of {G}ibrat's law},
	volume = {461},
	year = {2016}}

@book{GR,
	author = {Gradshteyn, I.S. and Ryzhik, I.M.},
	publisher = {Academic press},
	title = {Table of integrals, series, and products},
	year = {2014}}

@article{capasso1978,
	author = {Capasso, V. and Serio, G.},
	doi = {10.1016/0025-5564(78)90006-8},
	journal = {Math. Biosci.},
	number = {1-2},
	pages = {43--61},
	publisher = {Elsevier},
	title = {A generalization of the {K}ermack-{M}c{K}endrick deterministic epidemic model},
	volume = {42},
	year = {1978}}

@article{liu1986,
	author = {Liu, W. and Levin, S.A. and Iwasa, Y.},
	doi = {10.1007/BF00276956},
	journal = {J. Math. Biol.},
	pages = {187--204},
	publisher = {Springer},
	title = {Influence of nonlinear incidence rates upon the behavior of {SIRS} epidemiological models},
	volume = {23},
	year = {1986}}

@article{gualandi2019,
	author = {Gualandi, S. and Toscani, G.},
	doi = {10.1016/j.physa.2019.04.260},
	journal = {Phys. A},
	pages = {221--234},
	publisher = {Elsevier},
	title = {Size distribution of cities: A kinetic explanation},
	volume = {524},
	year = {2019}}

@article{szekely2013,
	author = {Sz{\'e}kely, G.J. and Rizzo, M.L.},
	doi = {10.1016/j.jspi.2013.03.018},
	journal = {J. Stat. Plan. Inference},
	number = {8},
	pages = {1249--1272},
	publisher = {Elsevier},
	title = {Energy statistics: A class of statistics based on distances},
	volume = {143},
	year = {2013}}

@article{bertaglia2023,
	author = {Bertaglia, G. and Pareschi, L. and Toscani, G.},
	doi = {10.3934/mbe.2024187},
	journal = {Math. Biosci. Eng.},
	number = {3},
	pages = {4241-4268},
	title = {Modelling contagious viral dynamics: a kinetic approach based on mutual utility},
	volume = {21},
	year = {2024}}

		\end{document}